\documentclass[notitlepage,leqno,10pt]{article}
\usepackage{latexsym}
\usepackage{amsmath}
\usepackage{amsfonts}
\usepackage{amssymb}
\usepackage{amscd}
\renewcommand{\a }{\alpha }

\renewcommand{\d}{\delta }
\newcommand{\D }{\Delta }

\newcommand{\e }{\varepsilon }
\newcommand{\g }{\gamma}
\renewcommand{\i }{\iota}
\newcommand{\G }{\Gamma }
\renewcommand{\l }{\lambda }
\renewcommand{\L }{\Lambda }

\newcommand{\n }{\nabla }

\newcommand{\Sig }{\Sigma}

\newcommand{\ov}{\overline}
\newcommand{\intbar}{\mathop{\int\makebox(-13.5,0){\rule[4pt]{.7em}{0.3pt}}%
\kern-6pt}\nolimits}
\newcommand{\wtilde }{\widetilde}

\newcommand{\be}{\begin{equation}}
\newcommand{\ee}{\end{equation}}
\newcommand{\bes}{\begin{equation*}}
\newcommand{\ees}{\end{equation*}}
\newcommand{\ba}{\begin{eqnarray}}
\newcommand{\ea}{\end{eqnarray}}
\newcommand{\bas}{\begin{eqnarray*}}
\newcommand{\eas}{\end{eqnarray*}}
\newenvironment{pf}{\noindent{\sc Proof}.\enspace}{\rule{2mm}{2mm}\medskip}
\newenvironment{pfn}{\noindent{\sc \bf Proof }}{\rule{2mm}{2mm}\medskip}

\newcommand{\R}{\mathbb{R}}

\newcommand{\N}{\mathbb{N}}

\author{ Mohameden  Ahmedou \; $\&$\;\; Mohamed Ben Ayed }

\date{}

\title{\bf  The Nirenberg problem on   high dimensional half spheres:\\
The effect of   pinching conditions}
\begin{document}

\newtheorem{lem}{Lemma}[section]
\newtheorem{pro}[lem]{Proposition}
\newtheorem{thm}[lem]{Theorem}
\newtheorem{rem}[lem]{Remark}
\newtheorem{cor}[lem]{Corollary}
\newtheorem{df}[lem]{Definition}

\maketitle

\bigskip

\centerline{   \emph{dedicated to the  memory of Prof. Louis Nirenberg }}

\begin{center}
{\bf Abstract}
\end{center}
In this paper we study the Nirenberg problem on standard half spheres $(\mathbb{S}^n_+,g), \, n \geq 5$, which  consists  of finding conformal metrics of prescribed scalar curvature  and zero boundary mean curvature on the boundary.  This problem amounts to solve the following boundary value problem involving the critical Sobolev exponent:
\begin{equation*}
(\mathcal{P}) \quad
\begin{cases}
   -\D_{g} u \, + \, \frac{n(n-2)}{4} u \,  =  K \, u^{\frac{n+2}{n-2}},\, u > 0 & \mbox{in } \mathbb{S}^n_+, \\
  \frac{\partial u}{\partial \nu }\, =\, 0  & \mbox{on } \partial \mathbb{S}^n_+.
\end{cases}
\end{equation*}
 where $K \in C^3(\mathbb{S}^n_+)$ is  a positive  function.\\
 This problem has  a variational structure but the  related Euler-Lagrange functional $J_K$ lacks compactness. Indeed it admits \emph{critical points at infinity}, which are \emph{limits} of  non compact orbits of the (negative) gradient flow. Through the construction of an appropriate  \emph{pseudogradient}  in the \emph{neighborhood at infinity}, we characterize these \emph{critical points at infinity}, associate to them an index, perform a \emph{Morse type reduction} of the functional $J_K$ in their neighborhood  and compute their contribution to the difference of topology between the level sets of  $J_K$, hence extending the  full Morse theoretical approach to this \emph{non compact variational problem}. Such an approach  is used to prove,  under various pinching conditions,   some existence results for  $(\mathcal{P})$ on half spheres of dimension $n \geq 5$.

\begin{center} \bigskip
\noindent{\bf Key Words:}  Non compact variational problems, Critical point at infinity, pseudogradient, Infinite dimensional Morse theory

\centerline{\bf AMS subject classification:  35C60, 58J60,  53C21.}

\end{center}

\tableofcontents

\section{Introduction and statement of the results}

In the early seventieth of the last century Louis Nirenberg asked the following question: Can a smooth positive function $K \in C^{\infty}(\mathbb{S}^n)$ defined on the standard $n-$dimensional sphere  $(\mathbb{S}^n, g)$ be realized as the scalar curvature of a metric $\ov g$ conformally equivalent to $g$ ?\\
On $\mathbb{S}^2$, setting $\ov g=e^{2u} g$ the Nirenberg problem is equivalent to solving the following nonlinear elliptic equation
\begin{equation*}
  -\D_{g} u \, + \, 1 \, = \, Ke^{2u}, \, \mbox{ in }\,  \mathbb{S}^2,
\end{equation*}
where $\D_{g}$ denotes the Laplace Beltrami operator.\\
For spheres of dimensions $n \geq 3$ and writing the conformal metric  as $\ov g:= u^{{4}/{(n-2)}} g$, the Nirenberg problem  amounts to solve the following nonlinear elliptic equation involving the Sobolev critical exponent:
\begin{equation}\label{eq:NP}
(\mathcal{NP}) \qquad -\D_{g} u \, + \frac{n(n-2)}{4} u \,  =  \, K u^{ \frac{n+2}{n-2}}\, ; \quad u > 0, \, \, \mbox{ in } \mathbb{S}^n.
 \end{equation}
The Nirenberg problem  has attracted a lot attention in the last half century. See \cite{KW1, Au76, BN83, AH91, BC1, CY, CGY, yyli1, yyli2, SZ, BCCH, Bahri-Invariant, BCH, CL1, CL2, Chen-Xu} and the references therein.
Actually
due to Kazdan-Warner obstructions, see \cite{KW1, BEZ1}, a positive answer to the Nirenberg's question requires imposing conditions on the function $K$. It turns out that finding  sufficient conditions under which the Nirenberg problem is solvable depends strongly on the dimension $n$ and the behavior of the function $K$ near its critical points. Indeed in low dimension $n< 5$ index counting criteria have been obtained, see \cite{BC1, CGY, Han, yyli1, yyli2}. Such a counting index criterium   fails,  under the nondegeneracy assumption  $(ND)$ (that is $\D K \neq 0$ at critical points of $K$), if the dimension $ n \geq 5 $. They can be extended on  high dimensional spheres in the  perturbative setting (that is  when $K$  is close to a constant) see \cite{CY, Chen-Xu} or under some flatness assumptions see \cite{yyli1, Ben, CL1}. \\
To explain the main difficulty in studying the Nirenberg problem and the differences between the low dimensional case $n< 5$ and the high dimensional one $n \geq 5$, we point out that due to the presence of the Sobolev critical exponent, the corresponding Euler-Lagrange functional does not satisfy the Palais-Smale condition. One way   to  overcome such a difficulty is to consider the following subcritical  approximation of the problem $(\mathcal{NP})$:
   \begin{equation}\label{eq:ep}
  (\mathcal{NP}_{\e}) \quad  -\D_{g} u \, + \, \frac{n(n-2)}{4} u \,  =  K \, u^{\frac{n+2}{n-2} - \e}, \quad  \, u > 0 \, \mbox{ in }\,  \mathbb{S}^n,
\end{equation}
where $\e > 0$ is a small parameter.
In this way one recovers the compactness and one then studies the behavior of  blowing up  solution $u_{\e}$ of $(\mathcal{NP}_{\e})$ as the parameter $\e$ goes to zero. Actually it can be proved that finite energy blowing up solutions of $ (\mathcal{NP}_{\e})$ can have only \emph{isolated simple blow up points} which are critical points of the function $K$, see \cite{yyli1, yyli2, CL2, MM19}. The reason of the  additional difficulty in the high dimensional case lies in the complexity of the blow up phenomenon. Indeed in dimensions $n=2,3$ there are only  single blow up points, see, \cite{CGY,Han,  BC1, yyli1, SZ} and  in dimension $n=4$ multiple bubbling may   occur only under some extra condition, see \cite{yyli2, BCCH} while, under the non degeneracy assumption $(ND)$, on spheres of dimension $n \geq 5$ every $m-$tuple $(q_1, \cdots, q_m)$ of distinct critical points  of $K$, satisfying $\D K(q_i)  < 0 $ for each $i =1, \cdots,m$  can be realized as a concentration set of blowing up solutions of  $(\mathcal{NP}_\e)$. See \cite{MM}.\\
Regarding the high dimensional case $n \geq 5$,  A. Malchiodi and M. Mayer \cite{MM19} obtained recently  an interesting  existence criterium under some  pinching condition. Their result reads as follows:\\
 {\bf Theorem A \cite{MM19}}
 {\it Let $n \geq 5$ and $K \in C^{\infty}(\mathbb{S}^n)$  be a  positive Morse function satisfying the following conditions
   $$  \forall q \in \mathbb{S}^n, \quad  \n K(q) = 0 \Rightarrow  \D K(q) \neq 0, \leqno{(i)} $$
   $$  {K_{max}}/{K_{min} } \, \leq  ({3}/{2})^{ {1}/{(n-2)}}, \,  \leqno{(ii)}$$
   where $K_{max} := \max_{\mathbb{S}^n} K$ and $K_{min} := \min_{\mathbb{S}^n} K$
$$ \# \{ q \in \mathbb{S}^n; \n K(q) = 0; \D K(q) < 0 \} \geq 2,  \leqno{(iii)} $$
where $\# A$ denotes the cardinal of the set $A.$ \\
 Then Nirenberg Problem $(\mathcal{NP})$ has at least one solution.}
\bigskip

In this paper we consider  a version of the Nirenberg problem on standard half spheres $(\mathbb{S}^n_+,g)$. Namely we prescribe simultaneously   the scalar curvature to be a positive function $ 0 < K \in C^3(\mathbb{S}^n_+)$ and the boundary mean curvature to be zero. This amounts to solve the following boundary value problem
\begin{equation}
(\mathcal{P}) \quad
\begin{cases}
   -\D_{g} u \, + \, \frac{n(n-2)}{4} u \,  =  K \, u^{{(n+2)}/{(n-2)}},\,  u > 0 & \mbox{ in } \mathbb{S}^n_+, \\
  \frac{\partial u}{\partial \nu }\, =\, 0  & \mbox{on } \partial \mathbb{S}^n_+,
\end{cases}
\end{equation}
 where $K \in C^3(\mathbb{S}^n_+)$ is a positive  function.\\
 This problem has been studied  on half spheres of dimensions $n=2, 3,4$. See  the papers \cite{LL, yyli, DMOA, BEOA, BEO, BOA, BGO} and the references therein.
Very much like the case of spheres,  to recover compactness  one considers  here  the  following subcritical approximation
 \begin{equation}
(\mathcal{P}_{\e}) \quad
\begin{cases}
   -\D_{g} u \, + \, \frac{n(n-2)}{4} u \,  =  K \, u^{\frac{n+2}{n-2} - \e}, \,  u > 0 & \mbox{ in } \mathbb{S}^n_+, \\
  \frac{\partial u}{\partial \nu }\, =\, 0  & \mbox{on } \partial \mathbb{S}^n_+.
\end{cases}
\end{equation}
Just as above, there are two alternatives for the behavior of a sequence of solutions $u_{\e}$ of $(\mathcal{P}_{\e})$. Either the $||u_{\e}||_{L^{\infty}}$ remains uniformly bounded  or it blows up and if it does $u^{{2n}/{(n-2)}}_{\e}  \mathcal{L}^n$ (where $ \mathcal{L}^n $ denotes the Lebesgue measure) converges to a sum of Dirac masses, some of them are sitting in the interior and the others ones are located on the boundary. The interior points are  critical points of $K$  satisfying  that $\D K \leq 0$ and the boundary points are critical points of $K_1$ the restriction of $K$ on the boundary and satisfying that ${\partial K}/{\partial \nu} \geq  0$. See \cite{DMOA, BOA, BEO}.
Furthermore  a refined blow up analysis, under the  non degeneracy assumption that $\D K \neq 0$ at interior critical points of $K$ and that ${\partial K}/{\partial \nu} \neq 0$ at critical points of $K_1$,  shows that in the dimension $n=3$  multiple bubbling may occur but all  blow up points are isolated simple, see \cite{yyli, DMOA}. Moreover in dimensions $n=2,3$ counting index criteria have been established, see \cite{LL, DMOA, BEO, BOA}. Furthermore   under additional condition on $K_1$ it has been proved in \cite{BGO} that all blow up  points are isolated simple, but already in dimension $n=4$  counting index formulae, under the above non degeneracy conditions fail. More surprisingly and in contrast with  the case of closed spheres, the Nirenberg problem on half spheres may have \emph{non simple blow up points}, even for finite energy bubbling solutions of  $(\mathcal{P}_\e)$ see \cite{AB}.

In this paper we study Problem $(\mathcal{P})$ from the viewpoint of \emph{the theory of critical points at infinity}. In this approach  initiated by  the late A. Bahri, see \cite{B1, BC1, BC2, Bahri-Invariant},  one studies the possible \emph{ends}  of non compact orbits of the (negative) gradient of  the associated Euler Lagrange functional. The method consists of taking advantage of the concentration-compactness analysis of non converging Palais-Smale sequences  to identify a potential \emph{neighborhood at infinity} where concentration may occur. Then one constructs  a  global \emph{ pseudogradient }  for which the  full analysis of the $\omega-$limit set,   in this neighborhood is easier than for the genuine gradient flow  and then uses it to characterize  \emph{critical points at infinity}. One then  performs a Morse reduction near these \emph{critical points at infinity} in order  to compute their topological contribution to the difference of topology between the level sets of the  Euler-Lagrange functional.\\
Before stating our main results, we set up some  notation and introduce our assumptions.\\
For the function $K$ and its restriction on the boundary  $K_1:= K_{\lfloor \partial \mathbb{S}^n}$, we use the following assumption: \\
$\textbf{(H1)}$: We assume that   $K$ is a  $C^3(\ov{\mathbb{S}^n_+} ) $  positive function, which  has only non-degenerate critical points with $\D K \neq 0$. (We point out that some of these points can be on the boundary.)\\
$\textbf{(H2)}$: We assume that the  restriction of $K$ on the boundary  $K_1:= K_{| \partial \mathbb{S}_+^n}$ has only non-degenerate critical points $z$'s. Furthermore we assume that  if $z$ is not a local maximum point of $K_1$, we have  that $\partial K/\partial \nu (z) \leq 0$.\\
$\textbf{(H3)}$:  If  $z \in \partial \mathbb{S}^n_+$ is a critical point  of $K_1$ satisfying that  $\partial K/\partial \nu (z) = 0$, hence $z$ is actually a critical point  of $K$   on $\partial \mathbb{S}^n_+$, we assume  that $\D K(z) \neq 0$ and  one of the following conditions is satisfied:
\begin{enumerate}
\item[(i)] either  $\partial K/\partial \nu (a) \D K(z) \leq 0$ for each $a\in \partial \mathbb{S}^n_+$ in a small neighborhood of  $z$,
\item[(ii)] or  $\lim_{a\in \partial \mathbb{S}^n_+;  a\to z} \frac{ \partial K/\partial \nu (a) }{d( a , z ) } =0$.\end{enumerate}
Next we introduce the following subsets of critical points of $K$ and $K_1$
\begin{align*}\label{setK}
&  \mathcal{K}^{-}_{in}    := \{ y\in \mathbb{S}^n_+: \n K (y)=0  \mbox{ and } \D K(y)<0\},  \\
 & \mathcal{K}^+_b  := \{ z\in \partial \mathbb{S}^n_+: \n K_1 (z)=0 \mbox{ and } \partial K/ \partial \nu (z ) > 0\} , \\
&  \mathcal{K}^{0,-}_b  := \{ z\in \partial \mathbb{S}^n_+: \n K_1 (z)=0 \, ; \,  \partial K/ \partial \nu (z ) = 0\mbox{ and } \D K(z)<0\} .
\end{align*}
Furthermore we define
$$
 \mathcal{K}^{\infty} \, := \,  \mathcal{K}^{-}_{in} \cup \mathcal{K}^+_b \cup \mathcal{K}^{0,-}_b .
$$

Our  first result is an existence result under a  pinching assumption, which parallels the above mentioned existence result of Malchiodi-Mayer. Namely we prove

\begin{thm}\label{t:pinching1}
Let $n \geq 5$ and  $0 < K \in C^3(\ov{\mathbb{S}^n_+})$ satisfying the assumptions  $(H1)$, $(H2)$ and $(H3)$.\\
If the following conditions hold
  $$ K_{\max} / K_{\min} < (5/4)^{1/(n-2)},\leqno{(i)} $$
  where $ K_{max} := \max_{\mathbb{S}^n_+} K$ and $K_{min} := \min_{\mathbb{S}^n_+} K $.
  $$  \# \mathcal{K}^{\infty} \geq 2, \leqno{(ii)} $$
  where $\# A$ denotes the cardinal of the set $A$.
Then Problem $(\mathcal{P})$ has at least one solution.
\end{thm}
\begin{rem}
\begin{enumerate}
  \item
  The above theorem is the counterpart of the existence result of Malchiodi-Mayer \cite{MM19}(see Theorem A quoted above). We point that the proof  of Theorem \ref{t:pinching1}, compared with the proof of Theorem A is more involved. In particular  the counting index argument in our case   is more subtle. Indeed due to the influence of the boundary the blow up picture is more complicated. Namely we have   boundary  and  interior blow up as well as mixed configurations involving both of them. Such a complicated picture  imposes to consider 4 critical levels instead of two critical levels needed in the case of  closed spheres. Such a fact  makes the index counting of the associated \emph{critical points at infinity} more involved, see Lemmas \ref{indicesy}, \ref{indicesz} in the appendix.
  \item
  The conditions $ (H2), (H3)$ are  used to rule out \emph{ non simple blow up}, see \cite{AB}. A phenomenon which does not occur in the case of closed spheres. See subsection \ref{s:simpleblow1}.
  \end{enumerate}
\end{rem}

\noindent
The above pinching condition $(i)$ of Theorem \ref{t:pinching1} can be relaxed when combined with some counting index formula involving either the boundary blow up points or the interior blow points. In the next theorem we provide an existence result involving the boundary blow up points. Namely we prove:

\begin{thm}\label{t:pinching2}
Let $n \geq 5$ and  $0 < K \in C^3(\ov{\mathbb{S}^n_+})$. Assume that the critical points of $K_1:= K_{|\partial \mathbb{S}^n_+}$ are non degenerate and that $K$ satisfies the assumption  $(H3)$.
If the following conditions hold
  $$ K_{\max} / K_{\min} < 2^{1/(n-2)} , \leqno{(a)} $$
  $$  A_1:=\sum_{z \in \mathcal{K}^+_b\cup \mathcal{K}^{0,-}_b} (-1)^{n-1 - morse(K_1,z)} \neq 1. \leqno{(b)} $$
Then Problem $(\mathcal{P})$ has at least one solution.
\end{thm}

Next we assume that  the above index formula $A_1 =1$, which implies, in particular  that the number of boundary blow up points is an odd number, say $2k+1$, where $k \in \N_0$. \\
The next existence result combined a pinching condition with a counting index formulae involving  interior blow up points. Namely we prove:

\begin{thm}\label{t:pinching3}
Let $n \geq 5$ and  $0 < K \in C^3(\ov{\mathbb{S}^n_+})$ satisfying the assumptions  $(H1)$, $(H2)$ and $(H3)$.\\
If the following conditions hold
  $$ K_{\max} / K_{\min} < (3/2)^{1/(n-2)} \mbox{ and } A_1 = 1,  \leqno{(i)} $$
  where $A_1$ is defined in Theorem \ref{t:pinching2},
  $$  B_1:= \sum_{y \in \mathcal{K}^-_{in}} (-1)^{n - morse(K,y)} \neq -k,  \leqno{(ii)} $$
  where
  $   \# (\mathcal{K}^{+}_b \cup \mathcal{K}^{0,-}_b) = 2k +1, \, k \in \N_0$. Then Problem $(\mathcal{P})$ has at least one solution.
\end{thm}

\bigskip
Regarding the method of proof of our main existence results, Theorems \ref{t:pinching1}, \ref{t:pinching2} and  \ref{t:pinching3} some comments are in order. Indeed although the general scheme falls in the framework of the techniques and ideas of  \emph{the critical point theory at infinity }, see \cite{Bahri-Invariant, BC1, BCCH}, the main arguments here are of a different flavor. Indeed with respect to the case of closed spheres, treated by A.Bahri   in  his seminal paper \cite{Bahri-Invariant}, the case of half spheres presents new aspects: From one part the blow up picture is more complicated (interior, boundary and \emph{mixed  configurations}) and from another part the behavior of the \emph{self interactions} of  interior bubbles and boundary bubbles is drastically  different. A fact which was used in \cite{AB} to construct subcritical solutions  having non simple blow ups. To rule out such a possibility, under our assumption $(H2)$ and $(H3)$,  we had to come up with a \emph{barycentric vector field } which moves a \emph{ cluster} of concentration points towards  their common barycenter  and to prove that along the flow lines of such a vector field the functional decreases and the concentration rates of an initial value do not increase, see Lemma \ref{simple}. Furthermore we prove that  in the neighborhood of \emph{critical points at infinity}, the concentration rates are comparable and the concentration points are not to close to each other. See  subsections \ref{s:simpleblow1} and \ref{s:simpleblow2}.



\bigskip

The remainder of this paper is organized as follows: In Section 2  we set up the variational framework and define the neighborhood at infinity and in Section 3 we construct an appropriate pseudogradient in the vicinity of highly concentrated bubbles and derive from the analysis of the behavior of its flow lines the set of its \emph{ critical points at infinity}. Section 4 is devoted to the proof of the main existence results of this paper. Lastly we collect in the appendix some  estimates of the bubble, fine asymptotic expansion of the Euler-Lagrange functional and its  gradient in the neighborhood at infinity as well as useful counting index formula for the critical points of the function $K$ and its restriction $K_1$ on the boundary.

\section{ Loss  of compactness and neighborhood at infinity}

In this section we  set up the analytical framework of the variational problem associated to the Nirenberg problem and recall the description of its lack of compactness.
Let $H^1(\mathbb{S}^n_+)$ be the Sobolev space endowed with the norm
$$
||u||^2 := \int_{\mathbb{S}^n_+} |\n u|^2 \, + \, \frac{n(n-2)}{4} \int_{\mathbb{S}^n_+} u^2,
$$
and let $\Sigma$ denote its  unit sphere.

Problem $(\mathcal{P})$ has a variational structure. Namely  its solutions are in one to one correspondence with the critical points of the functional
$$
J_K(u):= \frac{||u||^2}{( \int_{\mathbb{S}^n_+}  K |u|^{{2n}/{(n-2)}} )^{ {(n-2)}/{n}}} \quad \mbox{ defined on   } \Sigma^+:= \{  u \in \Sigma; \, u \geq 0 \}.$$
The functional $J_K$ fails to satisfy the Palais Smale condition. To describe non converging Palais-Smale sequences we introduce the following notation.\\
For $a \in \ov{ \mathbb{S}^n_+}$ and $\l > 0$ we define the \emph{standard bubble} to be
$$
\d_{a,\l}(x) := c_0 \frac{\l^{n-2/2}}{ (  \l^2 + 1 + (1-\l^2) \cos d(a,x))^{n-2/2}},
$$
where $d$ is the geodesic distance on  $\mathbb{S}^n_+$ and $c_0$ is a constant chosen such that
$$
- \D \d_{a,\l} \, + \, \frac{n(n-2)}{4} \d_{a,\l} \, = \, \d_{a,\l}^{{(n+2)}/{(n-2)}} \quad \mbox{ in } \mathbb{S}^n_+.
$$
For $a \in \ov{ \mathbb{S}^n_+},$ we define \emph{projected bubble} $\varphi_{a,\l}$ to be the unique solution of
$$
- \D \varphi_{a,\l} \, + \, \frac{n(n-2)}{4} \varphi_{a,\l} \, = \, \d_{a,\l}^{{(n+2)}/{(n-2)}} \quad \mbox{ in }\,  \mathbb{S}^n_+; \quad \frac{\partial \varphi_{a,\l}}{\partial \nu} \, = 0 \mbox{ on } \partial \mathbb{S}^n_+.
$$
We point out that $\varphi_{a,\l} = \d_{a,\l}$ if $a \in  \partial \mathbb{S}^n_+$.\\
Next for $m \in \N$ and $p,q \in \N_0$ such that $q + 2p = m$ we define the  {\em neighborhood of potential critical points at Infinity} $V(m,q,p,\e)$ as follows:
\begin{align*}
{V}(m,q,p, \e) := & \Big \{  u \in \Sigma:  \,
 \exists \,  \l_1, \cdots, \l_{p+q} > {\e^{-1}}; \,  \exists \, a_1, \cdots,a_{q+p} \in \ov{ \mathbb{S}^n_+},  \, \mbox{ with }\l_i d(a_i, \partial \mathbb{S}^n_+) < \e,  \, \forall \,  i \leq q,  \mbox{and }\\
&  \l_i d(a_i,\partial \mathbb{S}^n_+) > \e^{-1} \, \,  \forall \,    i > q,  \quad   \e_{ij} < \e \mbox{ such that }
 \| u - \frac{\sum_{i = 1}^{p+q}  K(a_i)^{(2-n)/4}\varphi_{a_i,\l_i} } { \| \sum_{i = 1}^{p+q}  K(a_i)^{(2-n)/4}\varphi_{a_i,\l_i}\|} \|< \e  \Big \},
\end{align*}
where
$$
\e_{ij}:= \Big( \frac{\l_i}{\l_j} + \frac{\l_j}{\l_i} + 2 \l_i \l_j  (1-\cos(d (a_i,a_j)))\Big)^{2-n/2}.
$$
In the following we describe  non converging Palais-Smale sequences. Such a description, which is by now standard,  follows from  concentration-compactness arguments as  in \cite{Lions, Struwe} and reads as follows

\begin{pro}
Let $u_k \in \Sigma^+$ be a sequence such that $\n J_K(u_k) \to 0$ and $J_K(u_k)$ is bounded. If Problem $(\mathcal{P})$ does not have a solution, then there exist $m\in \N$ and $p,q \in \N$ with $q + 2p = m$, a sequence of  positive real numbers $\e_k \downarrow 0$  as well as  subsequence of $u_k$, still denoted $u_k$ such that $u_k \in V(m,q,p,\e_k)$.
\end{pro}
Following A. Bahri and J-M. Coron, we consider for  $u \in V(m,q,p,\e)$ the following minimization problem

\begin{equation}\label{eq:min}
  Min  \left   \{  \|  u  \, - \, \sum_{i=1}^{p+q} \a_i \varphi_{a_i,\l_i} \|; \a_i > 0, \l_i > 0, a_i \in \partial \mathbb{S}^n_+, \forall i =1, \dots, q; \, a_i \in \mathbb{S}^n_+, \forall q+1 \leq i \leq q+p  \right  \}.
\end{equation}
We then have the following proposition whose proof is identical, up to minor modification to  the one  of  Proposition 7 in \cite{BC2}

\begin{pro}\label{p:min}
For any $m \in \N$ there exists $\e_m > 0$ such that if $\e < \e_m$ and $u \in V(m, q, p, \e)$ the minimization problem \eqref{eq:min}  has, up to permutation, a unique solution.
\end{pro}
Hence it follows from Proposition \ref{p:min} that every $u \in V(m, q,p, \e)$ can be written in a unique way as
\be \label{***}
u \, = \, \sum_{i=1}^{q} \a_i \d_{a_i,\l_i} \, + \,  \sum_{i= q+1}^{p+q} \a_i \varphi_{a_i,\l_i} \, + \, v,
\ee
where
$$
a_i \in \partial \mathbb{S}^n_+, \, i =1, \cdots,q \mbox{ and } a_i \in \mathbb{S}^n_+, \, i=q+1, \cdots, p+q,
$$
and  $v \in H^1(\mathbb{S}^n_+)$ satisfying
\begin{equation}\label{eq:V0}
(V_0) \quad \| v\| <  \e, \quad <v, \psi> = 0, \mbox{ for } \psi \in  \bigcup_{ 1 \leq i\leq q; \, q + 1 \leq j \leq q+p } \{  \d_i, \frac{\partial \d_i}{\partial \l_i}, \frac{\partial \d_i}{\partial a_i},  \varphi_j, \frac{\partial \varphi_j}{\partial \l_j}, \frac{\partial \varphi_j}{\partial a_j} \},
\end{equation}
where $ \d_i := \d_{a_i,\l_i}$ and $\varphi_i:= \varphi_{a_i,\l_i}$. In addition, the variables $\a_i$'s satisfy
\be \label{alpha} | 1 - J(u)^{n/(n-2)} \a_i^{4/(n-2)} K(a_i) | = o_\e (1) \quad \mbox{ for each } i.\ee

In the next lemma we deal with the $v$-part of $u \in V(m,q,p,\e)$ in order to prove, that its effect is negligible with the concentration phenomenon. Namely we prove:
\begin{lem}\label{eovv} Let $n \geq 5$.  For $\e > 0$ small,  there exists a $C^1$-map which, to each $(\a:= (\a_1,\cdots,\a_{p+q}) , a:= (a_1,\cdots,a_{p+q}), \l:= (\l_1, \cdots,\l_{p+q}) )$,  such that $u \, =\,  \sum_{i=1}^{p+q} \a_i \varphi_i \in V(m,q,p,\e)$,   associates $\overline{v}=\overline{v}_{(\a ,a,\l )}$ satisfying
$$
J_K( \sum_{i=1}^{p+q} \a _i \varphi_{a_i,\l_i}  +\overline{v})= \min\{ J_K( \sum_{i=1}^{p+q} \a _i \varphi_{a_i,\l_i}  +v) , \, v \mbox{ satisfies } (V_0)\}.
$$
Moreover, there exists $c>0 $ such that the following holds
$$ \|\overline{v} \| \leq c  \sum_{i=1}^{q+p} \frac{| \n K(a_i) | }{\l _i} + \frac{1}{\l _i^2} + \begin{cases}
& \sum \e _{ij }^{\frac {n+2}{2(n-2)}}(\ln \e _{ij}^{-1})^{\frac{n+2}{2n}} + \sum_{i > q }\frac{\ln(\l_i d_i)}{(\l_i d_i)^{(n+2)/2}}\mbox{ if } n \geq 6,\\
& \sum \e _{ij }(\ln \e _{ij}^{-1})^{{3}/{5}} + \sum_{i > q} \frac{1}{(\l_i d_i)^{3}}\mbox{ if } n = 5. \end{cases}
$$
\end{lem}

\begin{pf}
The proof follows as  in Proposition 3.1 in \cite{BCH}.  Indeed, easy computations imply that
\begin{align*} & J_K(u+v) = J_K(u) - f(v) +(1/2) Q(v) + o( \| v \|^2 ) \qquad \mbox{ where }\\
& f(v):= \int_{\mathbb{S}_+^n} K u^{\frac{n+2}{n-2}} v \quad \mbox{ and } \quad Q(v):= \| v \|^2 - \frac{n+2}{n-2} \sum_{i=1}^N \int_{\mathbb{S}_+^n}  \d_i ^{4/(n-2)} v^2.\end{align*}
Note that $Q$ is a positive definite quadratic form (see \cite{B1}) and we have that
\be\label{777} f(v) = \sum \a_i^{\frac{n+2}{n-2}} \int_{\mathbb{S}_+^n} K \varphi_i ^{\frac{n+2}{n-2}} v + O \Big( \sum_ {i\neq j } \int_{\mathbb{S}_+^n} \sup(\varphi_j, \varphi_i) ^{\frac{4}{n-2}} \inf(\varphi_j, \varphi_i)  | v | \Big).\ee
Observe that, for $n\geq 6$, it follows that $4/(n-2) \leq 1$. Hence, using Holder's inequality, we get
\begin{align}
\int_{\mathbb{S}_+^n} \sup(\varphi_j, \varphi_i) ^{\frac{4}{n-2}}  & \inf(\varphi_j, \varphi_i)  | v |  \leq   \int_{\mathbb{S}_+^n} (\varphi_j \varphi_i) ^{\frac{n+2}{2(n-2)}} | v | \nonumber\\
& \leq c \| v \| \Big(  \int_{ \mathbb{S}_+^n} (\d_j \d_i) ^{\frac{n}{n-2}} \Big)^{\frac{n+2}{2n}} \leq c  \| v \| \e _{ij }^{\frac {n+2}{2(n-2)}}(\ln \e _{ij}^{-1})^{\frac{n+2}{2n}} \quad \mbox{ if } n \geq 6,\label{ww1}\\
\int_{\mathbb{S}_+^5} \sup(\varphi_j, \varphi_i) ^{{4}/{3}}  & \inf(\varphi_j, \varphi_i)  | v |  \leq c \| v \| \e _{ij }(\ln \e _{ij}^{-1})^{{3}/{5}} \mbox{ if } n = 5.\label{ww2}
 \end{align}
For the other term, for $i\leq q$ (that is $a_i \in \partial \mathbb{S}_+^n$), using the fact that $\langle \d_i, v \rangle =0$, we get
$$ \int_{\mathbb{S}_+^n} K \d_i ^{\frac{n+2}{n-2}} v =  O\Big( | \n K(a_i) | \int_ {\R^n_+} | x-a_i|  \d_i ^{\frac{n+2}{n-2}} | v | +  \int_ {\R^n_+} | x-a_i|^2  \d_i ^{\frac{n+2}{n-2}} | v |   \Big)   = O \Big(\Big( \frac{| \n K(a_i) | }{\l _i} + \frac{1}{\l _i^2} \Big)\| v\|\Big).$$
For $i > q$, using Lemma \ref{lem:1}, we get
\begin{align*}
  \int_{\mathbb{S}_+^n} K \varphi_i ^{\frac{n+2}{n-2}} v  & =   \int_{\mathbb{S}_+^n} K \d_i ^{\frac{n+2}{n-2}} v +  O\Big(\int_{\mathbb{S}_+^n} \d_i ^{\frac{4}{n-2}} | \varphi_i - \d_i | | v | \Big) \\
& = O \Big(\Big( \frac{| \n K(a_i) | }{\l _i} + \frac{1}{\l _i^2} \Big)\| v\|\Big) + \begin{cases}  O \big(  \| v \|  /(\l_i d_i ) ^{n-2}  \big)\mbox{ if } n = 5,\\
                                                                            O \big(  \| v \| \ln(\l_i d_i) /(\l_i d_i ) ^{(n+2)/2} \big) \mbox{ if } n \geq 6 \end{cases} \end{align*}
and the result follows.
\end{pf}



\section{ Pseudogradient and  Morse Lemma  at Infinity }

This section is devoted to the construction of a   pseudogradient for the functional $J_K$, which has the   property that along its flow lines there could be only finitely many isolated blow up ponits. Such a pseudogradient coincides with the gradient outside of $\bigcup_{m,q,p} V(m,q,p,\e/2)$ and satisfies   the Palais-Smale condition there.  Moreover in each $V(m,q,p,\e)$ it has the property to move the concentration points according to $\n K$ or $\n K_1$, the $\a_i$'s to their maximum values and the concentration $\l_i$'s are moved so that the functional $J_K$ decreases along its  flow lines. The global vector field is then defined by  convex combining these two vector fields.  Such a construction is then used to perform a Morse reduction near the singularities of the pseudogradient and to compute the difference of topology induced by the \emph{ critical points at infinity}  between the level sets of the Euler-Lagrange functional $J_K$.\\
The first step in the construction of the pseudogradient is to describe the movement of  the variable $v$. In fact,  since $\ov v$ minimizes $J_K$ in the $v$-space, it follows from the classical Morse Lemma that   there exists a change of variable $v \to V$ such that
\be \label{vV} J_K(\sum_{i=1}^{p+q} \a_i \varphi_{a_i,\l_i} + v ) = J_K(\sum_{i=1}^{p+q} \a_i \varphi_{a_i,\l_i} + \ov{v} ) + \| V \| ^2.\ee
Hence, for the variable $V$, we will use $\dot{V} = -V$ to bring it to $0$. Thus, we need to construct some vector fields by moving the variables $\a_i$, $a_i$ and $\l_i$.

\subsection{ The case of a  single concentration  point}

We point out that the construction of a pseudogradient satisfying  the above properties becomes quite involved in the case of more than one concentration point. Indeed in the case of two bubbles sitting at different points, their mutual interaction comes into play. For this reason we start by constructing the needed pseudogradient in  neighborhoods at infinity, containing one interior or one boundary point. To do so   we  consider two cases, the first one corresponds to  $p=1$ and $q=0$ (case of an interior concentration point) and the second one corresponds to $p=0$ and $q=1$ (the case of a boundary point). Namely we prove:

\begin{pro} \label{pp:champ2m1} Assume that $K$ satisfies $(H1)$ and $(H3)$.  A  pseudogradient $W$  can be defined  so that the following holds : There is a constant $c>0$ independent of $u=\a \varphi_{a,\l} \in V(2p+q, q,p,\e)$  (with $q=1$ or $p=1$) such that
\begin{align*}
(i)\,\,  &\langle-\nabla J_K(u),W\rangle \geq  c \begin{cases} & 1/\l^2 + 1/(\l d)^{n-2} + | \n K(a) | /\l \mbox{ if } p=1; q = 0,\\
                           & 1/\mu + |1-J(u)^{\frac{n}{n-2}}\a^{\frac{4}{n-2}} K(a)| \mbox{ if } p=0; q = 1, \end{cases} \\
 (ii)  \, \,   &\langle-\nabla J_K(u+\overline{v}),W+\frac{\partial\overline{v}}{\partial(\a,a,\l)}(W)\rangle \geq  c \begin{cases} & 1/\l ^2 + 1/(\l d)^{n-2} + | \n K(a) | /\l \mbox{ if } p=1; q = 0,\\
                           & 1/\mu + |1-J(u)^{\frac{n}{n-2}}\a^{\frac{4}{n-2}} K(a)| \mbox{ if } p=0; q = 1, \end{cases} \end{align*}
where $d:= d(a,\partial \mathbb{S}^n_+)$ for $a \in  \mathbb{S}^n_+$ and  $ \mu^{-1}= |\n K(a) | /\l + 1/\l^2 $ for $a\in \partial \mathbb{S}^n_+$.\\
$(iii)$ The vector field $W$ is bounded  with the property that along its flow lines,  $\l$  increases only in  the following region
\begin{itemize}
  \item
  If $p=1$ then  $\l$ increases  if and only if the point $a $ belongs to a  small  neighborhood  of  a critical point   $y\in \mathbb{S}_+^n$  of $K$, such that  $\D K(y) < 0$
  \item
   If $q=1$ then  $\l$ increases if and only if the point $a $  belongs to a small  neighborhood  of  a critical point  $z \in \partial \mathbb{S}_+^n$ of $K_1$ such that either  $(\partial K/\partial \nu) (z)>0$ or $(\partial K/\partial \nu) (z)=0$ and $\D K(z) < 0$.
\end{itemize}
 \end{pro}

\begin{pf}
We start by giving the proof of Claim $(i)$ for the case where $p=1$ and $q=0$ that is in $V(2,0,1,\e)$.
First, we notice  that, if $a$ is close to a critical point  $y$ of $K$ in $\mathbb{S}_+^n$, then $\D K(a) = \D K(y) (1+o(1))$ and therefore $\D K(a)$ has a constant sign. \\
 Let $M$ be a large constant and let  $\psi_1$ be a $C^\infty$ cut off  function defined by $\psi \in [0,1]$, $\psi_1(t) = 1 $ if $t \geq 2$ and $\psi_1(t)=0$ if $t\leq 1$. We define

$$ W := \psi_1\Big(\frac{ \l | \n K(a) | }{M} \Big) \Big( \frac{1}{\l} \frac{\partial \varphi_{a,\l}}{\partial a} \frac{ \n K(a)} {| \n K(a) |} - \l  \frac{\partial \varphi_{a,\l}}{\partial \l} \Big) + \Big(1 - \psi_1\Big(\frac{ \l | \n K(a) | }{M} \Big)\Big) (\mbox{ sign}(-\D K(a))) \l \frac{\partial \varphi_{a,\l}}{\partial \l} .$$
We notice that, in the region where $| \n K(a) | \geq 2M/\l$, we have  that $\psi_1( { \l | \n K(a) | }/{M}  )=1$, therefore the  Claim $(i)$ follows from Proposition \ref{pp65}. \\
Next if  $| \n K(a) | \leq 2M/\l$ then  $a$ is very  close to a critical point  of $K$ in $\ov{ \mathbb{S}^n_+}$. We claim that this critical point  cannot be on the boundary. Indeed, arguing by contradiction, we assume that  $a$ is in small neighborhood  of a critical point  $z\in \partial \mathbb{S}_+^n$.  Since  $z$ is a non-degenerate critical point  of $K$, we   derive that  $\l d(a,z)$ is bounded which contradicts the fact that $\l d(a,\partial \mathbb{S}_+^n)$ is very large.  Hence our claim follows and $a$ is close to an interior critical point  $y$ in $\mathbb{S}_+^n$. \\
Next  using Proposition \ref{pp65} we derive that
$$ \langle-\nabla J_K(u),W\rangle \geq  c  \psi_1\Big(\frac{ \l | \n K(a) | }{M}\Big ) \Big( \frac{ | \n K(a) | }{\l } + \frac{ 1 }{\l ^2}\Big) +  \Big(1 - \psi_1\Big(\frac{ \l | \n K(a) | }{M} \Big)\Big) \frac{ c }{\l^2}$$ which implies Claim $(i)$ in this region.\\
 Hence Claim $(i)$ is proved in the case where $p=1$ and $q =0$.\\
 Concerning $(ii)$  it follows from $(i)$ using the estimate of $\ov v$ in Lemma  \ref{eovv}.
 Finally we notice that    $\l$ increases along the flow lines of  the pseudogradient $W$ only in the region  where $a$ is close to a critical point  $y$ with $\D K(y) < 0.$ Thus the proof of the proposition follows  in the case where $p=1$ and $q=0$.

  Next we consider the case where $p=0$ and $q=1$, that is the case of a boundary  concentration point $a \in \partial \mathbb{S}_+^n$. In this situation we   divide the set  $V(1,1,0,\e)$ into 3 subsets and construct an appropriate vector field in each of these sets.\\
 $\textbf{(1)}$ Let $V_1^1:= \{ u \in V(1,1,0,\e) :    | 1- J_K(u)^{\frac{n}{n-2}}\a ^{\frac{4}{n-2}}K(a) | \geq   M /\mu   \}$. In this region, we define
 $$  \underline{W}^1_1:=  \mbox{ sign}( 1- J_K(u)^{\frac{n}{n-2}}\a ^{\frac{4}{n-2}}K(a) | )\, \,  \d_{a,\l}$$ and using Proposition \ref{p:35}, Claim $(i)$ follows easily (since $M$ is chosen large).\\
 $\textbf{(2)}$ Let $V_1^2:= \{ u \in V(1,1,0,\e) :    | 1- J_K(u)^{\frac{n}{n-2}}\a ^{\frac{4}{n-2}}K(a) | \leq   2M /\mu  \mbox{ and }  | \n K_1(a) | \geq \eta \}$,  where $\eta$ is a small fixed constant. In this region, we define
$$ \underline{W}^2_1:=  \frac{1}{\eta} W_a^{b} \quad \mbox{ where }  W_a^{b}:= \frac{1}{\l} \frac{\partial \d_{a,\l} }{\partial a} \frac{ \n K_1(a) }{| \n K_1(a) | } .$$
Note that, in this region, the parameter $\mu$ is of the same order that $\l$. Hence, using Proposition \ref{p:34}, the proof of Claim $(i)$ follows. \\
 $\textbf{(3)}$ Let $V_1^3:= \{ u \in V(1,1,0,\e) :    | 1- J_K(u)^{\frac{n}{n-2}}\a ^{\frac{4}{n-2}}K(a) | \leq   2M /\mu  \mbox{ and }  | \n K_1(a) | \leq 2 \eta \}$. In this region, $a$ is close to a critical point $z$ of $K_1$. The pseudogradient will depend on $z$. We define
 \be\label{undW13} \underline{W}^3_1:= \psi_1( {\l | \n K_1(a) | } / {M}) W_a^b + \gamma \, \, \l \frac{\partial \d _{a,\l}}{\partial \l}  \quad \mbox{ with } \g\in \{-1,1\} \mbox{ satisfying }\ee
$$  \begin{cases}
         & \gamma = 1 \quad \mbox{if } \partial K/\partial \nu (z) > 0 \mbox{ or } \partial K/\partial \nu (z) = 0 \mbox{ and  } \D K(z) < 0, \\
         & \gamma = - 1 \quad \mbox{if } \partial K/\partial \nu (z) < 0 \mbox{ or } \partial K/\partial \nu (z) = 0 \mbox{ and  } \D K(z) > 0.\end{cases}$$
Using Propositions \ref{p:33} and \ref{p:34}, it holds
\be\label{mm1}  \langle-\nabla J_K (u), \underline{W}^3_1\rangle \geq  c  \psi_1\Big(\frac{ \l | \n K_1(a) | }{M} \Big) \Big( \frac{ | \n K_1(a) | }{\l } + \frac{ 1 }{\l ^2}\Big) + \gamma \Big( \frac{c_3}{\l} \frac{\partial K}{\partial \nu}(a) - c \frac{\D K(a)}{\l^2} + O(\frac{1}{\l^3})\Big).\ee
Observe that, if $\partial K/\partial \nu (z) \neq 0$, it follows that  $\gamma \partial K/\partial \nu (a) \geq c > 0$ and therefore Claim $(i)$ follows easily. In the other case, that is $\partial K/\partial \nu (z) = 0$, we need to make use  of the assumption $(H3)$. Indeed,
\begin{itemize}
\item if $(i)$ of $(H3)$ holds, it follows that $\gamma \partial K/\partial \nu (a) = |  \partial K/\partial \nu (a) |$ and $ - \gamma \D K(a)\geq c > 0$. Therefore, if $\l |\n K_1(a) | \geq 2 M$, in the lower bound of  \eqref{mm1} will appear $| \n K_1(a) | /\l + |  \partial K/\partial \nu (a) | / \l + 1 /\l^2 $ which is larger than $c /\mu$. Hence, Claim $(i)$ follows in this case. However, if $\l |\n K_1(a) | \leq 2 M$, it follows that $| \n K(a) |   \leq c M /\l$ (since we assumed that $z$ is a non degenerate critical point). Therefore $1/\l^2 \geq c (1/\l^2 + |\n K(a)|/\l) = c/\mu$. Thus Claim $(i)$ follows in this case.
\item Next  we consider the case where $(ii)$ of $(H3)$ holds. Recall that $z$ is a non degenerate critical point  of $K_1$, thus it follows that there exists $r_1 > 0$ such that $| \n K_1(a) | \geq \underline{c} d(a,z)$ for each $a \in B(z,r_1)$. Let $\varrho_1 > 0$ (satisfying $\varrho_1 \max(M,1/\underline{c})$ is very small), using $(ii)$ of $(H3)$, there exists $r_2>0$ (with $r_2 \leq r_1$) such that $| \partial K/ \partial \nu (a) | \leq \varrho_1 d(a,z)$  for each $a\in B(z,r_2)$. Hence, in $B(z,r_2)$,  $| \partial K/ \partial \nu (a) | = o(| \n K_1(a) | )$ (since $\varrho_1$ is chosen so that $\varrho_1 / \underline{c}$ is small)  and therefore $| \n K_1(a) | = | \n K(a) | (1+o(1))$. Finally, as before, if $\l |\n K_1(a) | \geq 2 M$, in the lower bound of  \eqref{mm1} will appear $| \n K_1(a) | /\l$. Furthermore, we have $ - \gamma \D K(a)\geq c > 0$ and $| \partial K/ \partial \nu (a) | = o(| \n K_1(a) | )$ which imply the proof of Claim $(i)$ in this case. In the other case, which is $\l |\n K_1(a) | \leq 2 M$, it holds: $d(a,z) \leq c M/\l$ which implies that
$| \partial K/ \partial \nu (a) | \leq \varrho_1 d(a,z) \leq c \varrho_1 M/ \l^2 = o(1/\l^2)$ (by the chose of $\varrho_1$). Thus the proof of Claim $(i)$ follows from \eqref{mm1}.
\end{itemize}

Finally Claim $(ii)$   follows from Claim $(i)$ using the estimate of $\ov v$ in Lemma \ref{eovv} and Claim $(iii)$ follows immediately from the properties of the constructed vector field.
\end{pf}

We remark that the assumption $(H2)$ is not used in the construction of the pseudogradient in $V(1,1,0,\e)$.

\subsection{ The case of multiple concentration points}

\noindent
In  the next proposition we address the case where the set of the concentration points contains more  than one point. Before stating our result we define  for $i =1, \cdots, m$  the scalar quantity $\mu_i$  as follows
\be\label{mui}
\mu_i^{-1}= |\n K(a_i) | /\l_i + 1/\l_i^2 \mbox{ if } i\leq q\, \, ; \quad  \mu_i= \l_i^2 \mbox{ if } i\geq q+1.
\ee
The behavior of such a quantity along the flow lines of the constructed pseudogradient plays crucial role in identifying \emph{critical points at infinity}.

\begin{pro} \label{pp:champ2} Assume that $K$ satisfies $(H1)$, $(H2)$ and $(H3)$.  A  pseudogradient $W$  can be defined  so that the following holds : There is a constant $c>0$ independent of
$u=\sum_{i=1}^q\a_i\d_{a_i,\l_i}+\sum_{j=q+1}^{p +q} \a_j \varphi_{a_j,\l_j}\in V(m, q,p,\e)$ such
that
\begin{align*}
(i)\qquad &\langle-\nabla J_K(u),W\rangle\geq  c \,  \sum_{i=1}^{p+q}  \frac{1}{\mu_i^{2-{1}/{(n-2)}}}+ c  \sum_{i \leq q} |1-J_K(u)^{\frac{n}{n-2}}\a_i^{\frac{4}{n-2}} K(a_i)|^{2-\frac{1}{n-2}}  \\
 &\qquad + c\, \sum_{k\neq r} \e_{kr}^{\frac{n-1}{n-2}}  +c \sum_{i> q}\Big(  \frac{1}{(\l_id_i)^{n-1}}+\Big(\frac{|\nabla K(a_i)|}{\l_i}\Big)^{2-\frac{1}{n-2}} \Big)\\
(ii) \qquad  &\langle-\nabla J_K(u+\overline{v}),W+\frac{\partial\overline{v}}{\partial(\a_i,a_i,\l_i)}(W)\rangle \geq  c\,  \sum_{i = 1}^{p+q} \frac{1}{\mu_i^{2-{1}/{(n-2)}}}+ c\, \sum_{k\neq r} \e_{kr}^{\frac{n-1}{n-2}}\\
& \qquad+c \sum_{i\leq q}   |1-J_K(u)^{\frac{n}{n-2}}\a_i^{\frac{4}{n-2}} K(a_i)|^{2-\frac{1}{n-2}}     +c \sum_{i> q}\Big( \frac{1}{(\l_id_i)^{n-1}}+\Big(\frac{|\nabla K(a_i)|}{\l_i}\Big)^{2-\frac{1}{n-2}} \Big)
\end{align*}
where $d_i:= d(a_i, \partial \mathbb{S}^n_+)$.\\
(iii) The vector field $W$ is bounded  with the property that along its flow lines  the maximum of the $\mu_i$'s  increases only if the $(q+p)-$tuple $(a_1,\cdots,a_q, \cdots a_{q+p})$ is close to a collection of different critical points of $K$ or $K_1$ $(z_1,\cdots,z_q,y_{q+1}, \cdots y_{q+p})$   with the $y_i$'s are critical points of $K$ in $\mathbb{S}^n_+$ satistying  $\D K(y_{i}) < 0$ for each $ i \geq q+1$ and the $z_i$'s  are critical points of  $K_1$ such that either  $(\partial K/\partial \nu) (z_{i_k})>0$ or $((\partial K/\partial \nu) (z_{i})=0$ and $\D K(z_{i}) < 0)$.
\end{pro}

The  construction of a pseudogradient satisfying  $(i)$,  $(ii)$,  $(iii)$  is quite involved and requires some preparatory Lemmas and estimates. Its construction depends on the behavior  of the leading terms of the $\a$-, $a$- and $\l$-component of the gradient in  the neighborhood at infinity $V(m,q,p,\e)$. To perform such a  construction we  divide the set $V(m,q,p,\e)$ into four  subsets. The first and the second ones correspond to the situation  where at least one of the variables $\a_i$'s and $a_i$'s is not in its  critical position and the $\mu_i$'s are of the same order. In the third one, the $\mu_i$'s are still of the same order but the variables $\a_i$'s and $a_i$'s are very close to their critical positions. Finally  in the fourth one we address the case where  the $\mu_i$'s are not of the same order. \\
To define these regions, we introduce the following notation. For  $M_2$  a large constant we set:
$$ \G_{\a_k}:= \frac{ | 1- J_K(u)^{\frac{n}{n-2}}\a  _k^{\frac{4}{n-2}}K(a_k) | }  {  M_2 (\sum_{r\neq k} \e_{kr} +1/\mu_k) } \, \, ; \, \, \G_{a_i}^b:= \frac{ | \n K_1 (a_i) | / \l_i }{  M_2 /\l_i^2 + (1/ M_2^2) \sum_{k\in I}  \e_{ik}}\quad \mbox{ for } i \leq q$$
\be \G_{a_i}:= \frac{ |\n K(a_i)| /\l_i  } {M_2 (\sum_{  k \neq i } \e_{ki} +{(\l_i d_i)^{2-n}}+ \frac{1}{\l_i^2})}\, ; \, \,  \G_{H_i}:= H(a_i, a_i) / M_2 \l_i^{n-4} \quad \mbox{ for } i > q ,\label{cq4}\ee
$$ \G_{\l_k} := \mu_k \sum_{j \neq k }\e_{jk} /M_2 \quad \mbox{ for } 1\leq i \leq q +p . $$
To explain the relevance of the above quantities, we state the following Lemma
\begin{lem}\label{pointint} 1) Let $a_i$ be an interior point satisfying   $\G_{\l_i} +\G_{a_i} +\G_{H_i}  \leq 8$. Then $a_i$ is close to a interior critical point  $y$ of $K$ in  $\mathbb{S}^n_+$.\\
2) If  $a_i, a_j$ are interior points satisfying that  $\G_{\l_k} +\G_{a_k} +\G_{H_k}  \leq 8$ for $k=i,j$ and if their  corresponding concentration rates  $\l_i$ and $\l_j$ are of the same order. Then  $a_i$ and $a_j$ cannot be close to the same critical point.
\end{lem}
\begin{pf} Since $i$ satisfies: $\G_{H_i} + \G_{a_i} + \G_{\l_i} \leq 8$, this implies that $| \n K(a_i) | \leq C/\l_i$ and therefore $a_i$ is close to a critical point  of $K$. We need to exclude the case where  this critical point  lies on the boundary. In fact, assuming that it is the case, i.e. $a_i$ is close to  $z\in \partial \mathbb{S}^n_+$. Then it follows from  $(H1)$, that  $\l_i d(a_i,z)$ is bounded, which is not allowed.  Therefore, each concentration point $a_i$ is close to a critical point  $y_{j_i} \in \mathbb{S}^n_+$ and  the first assertion is proved.\\
Concerning the second one, assume that two different points $a_i$ and $a_j$ are near the same critical point $y$. Then we have from the first assertion: $\l_k d(a_k,y)$ is bounded for $k=i,j$. Since  $\l_i$ and $\l_j$ are assumed to be of the same order, it follows that $\l_k d(a_i , a_j ) $ is bounded, which  contradicts  the smallness of  $\e_{ij}$.
\end{pf}

\subsubsection{Construction of some local pseudogradients}
In this subsection we construct some local pseudogradients in some parts of the neighborhood at infinity. These vector fields will be glued together to obtain a global pseudogradient satisfying the properties required in Proposition \ref{pp:champ2}.\\
For $M_0$ a large number we define   the following subsets of $V(m,q,p,\e)$
\begin{align*}
 V_1(M_0):= & \{ u : \mu_{\max}  \leq   2M_0 \,  \mu_{\min} \} \cap  \{ u : \exists \, \, i > q : \G_{H_i} + \G_{a_i} + \G_{\l_i} \geq 6 \} ,\\
 V_2(M_0):= & \{ u : \mu_{\max} \leq  2M_0\,  \mu_{\min} \}  \cap \{ u : \forall \, \, i > q  : \G_{H_i} + \G_{a_i} + \G_{\l_i} \leq 8 \}  \cap \Big( \{ u : \exists \, \, i \leq q  : \G_{\a_i}  + \G_{\l_i} \geq 4 \}\\
 &  \cup \{ u: \exists \, \, i \leq q : d(a_i ,\mathcal{K}^b) \geq \eta\}\Big)  \qquad \mbox{ where } \mathcal{K}^b:= \{ z\in \partial \mathbb{S}^n_+: \n K_1(z)=0\} ,\\
 V_3(M_0):= & \{ u : \mu_{\max} \leq  2M_0\,  \mu_{\min} \}  \cap \{ u : \forall \, \, i > q  : \G_{H_i} + \G_{a_i} + \G_{\l_i} \leq 8 \}  \cap \{ u : \forall \, \, i \leq q  : \G_{\a_i}  + \G_{\l_i} \leq 6 \}\\
 &  \cap \{ u: \forall \, \, i \leq q : d(a_i ,\mathcal{K}^b) \leq 2\eta\} , \\
 V_4(M_0):= & \{ u : \mu_{\max} > M_0\,  \mu_{\min} \}, \end{align*}
 where
 $ \mu_{\max}:= \max_j \mu_j$   and  $ \mu_{\min}:= \min_j \mu_j$.

Before  defining a pseudogradient in each subset, we single out some of their properties that will be used in the construction of the local pseudogradients.
\begin{rem}\label{ree1}
1) In $V_k(M_0)$, for $k\leq 3$, the variables $\mu_i$'s are of the same order. Thus, using Lemma \ref{epsij1}, we derive that, for each $i\neq j \leq q$, it holds
\be\label{bb1} - \l_i \frac{\partial \e_{ij}}{\partial \l_i} \geq c \e_{ij} .\ee
Furthermore, for $i\ne j > q$, we deduce that $\l_i$ and $\l_j$ are of the same order and therefore \eqref{bb1} holds true. Now, for $i > q$ and $j \leq q$, we have $\l_i d_i$ is very large which implies that $\l_i d(a_i ,a_j )$ is also very large and therefore \eqref{bb1} holds for these indices.

 2) In $V_k(M_0)$, $k=2,3$, for each $i > q$, the concentration point $a_i$ is close to a critical point  $y_{j_i} \in \mathbb{S}^n_+$ and two different points $a_i$ and $a_j$ cannot be near the same critical point $y$ (see Lemma \ref{pointint}).

 3) In $V_3(M_0)$,  for each $i \leq q$,  $a_i$ is close to a critical point  $z_{j_i}$ of $K_1$ in $\partial \mathbb{S}^n_+$.
 \end{rem}
We start our construction by defining a pseudogradient in $V_1(M_0)$.
\begin{lem}\label{chv1}
There exists a  bounded pseudogradient $W_1$  so that the following holds : There is a constant $c>0$ independent of
$u=\sum_{i=1}^{q}\a_i\d_i  + \sum_{i=q+1}^{p+q}\a_i\varphi_i \in V_1(M_0)$ such that
\be\label{estchv1} \langle-\nabla J_K(u),W_1\rangle \geq   \,  \sum_{i=1}^{q+p} \frac{c}{\mu_i}+ c \sum_{i=1}^{q}|1-J_K(u)^{\frac{n}{n-2}}\a_i^{\frac{4}{n-2}} K(a_i)|  + c\, \sum_{k\neq r} \e_{kr}  + \sum_{i=q+1}^{q+p} \frac{| \n K(a_i)|}{\l_i} . \ee
Furthermore, the $\l_i$'s are decreasing functions along the flow lines generated by this pseudogradient. In addition, the constant of $1/\mu_{\max}$ is independent of $M_0$ and $M_2$.
\end{lem}
\begin{pf}  
We start by defining the following vector fields:
\be \label{aaz1} W_{\L _{in}} := - \sum_{i > q} (\psi_1(\G_{\l_i}) +\psi_1(\G_{H_i})) \l_i \frac{\partial \varphi_i}{\partial \l_i} \quad \mbox{ and } \quad W_a^{in}:=  \sum _{i > q }   \psi_1(\G_{a_i})  \frac{1}{\l_i}  \frac{\partial \varphi_i}{\partial a_i} \frac{\n K(a_i)}{| \n K(a_i)|} \ee
\be \label{aaz2} W_{\L _{b}} := - \sum_{i \leq  q} \psi_1(\G_{\l_i}) \l_i \frac{\partial \d_i}{\partial \l_i} \quad \mbox{ and } \quad {W}_\alpha:= - \sum_{k \leq q }  \psi_1 (\G_{\a_k})\mbox{ sign}(1- J_K(u)^{\frac{n}{n-2}}\a  _k^{\frac{4}{n-2}}K(a_k)) \d_k \ee
where $\psi_1$ is a $C^\infty$ function defined by $\psi_1 \in [0,1]$, $\psi_1(t) = 1 $ if $t \geq 2$ and $\psi_1(t)=0$ if $t\leq 1$.
 Observe that, using Propositions  \ref{p:33}, \ref{pp65}, the estimate  \eqref{bb1} and the definition of $\psi_1$, we derive that
\begin{align}
& \langle - \n J_K(u), W_{\L_{in}} \rangle \geq c   \sum_{ i > q }(\psi_1(\G_{\l_i})+\psi_1(\G_{H_i})) \Big( \sum_{j\neq i} \e_{ij}  +  \frac{H(a_i,a_i)}{\l_i^{n-2}} + \frac{M_2}{2}\frac{1}{\l_i^2} + O(R_1)  \Big) := \ov{ \G}_{\L_{in}} \,   \label{cq22}\\
& \langle -\n J_K(u), W_{\L_{b}} \rangle  \geq c \sum_{i \leq q}\psi_1(\G_{\l_i}) \Big( \sum_{ k\neq i} \e_{ik} +\frac{M_2}{2} \frac{1}{\mu_{i}} +O\Big(  \sum_{ k > q } \e_{ki} + R_1^b \Big) \Big) := \ov{\G}_{\L_{b}} .\label{cq222}
\end{align}
Moreover   using  Proposition \ref{p:35}, we derive that
\be \label{estchalpha}
\langle -\n J_K(u), {W}_\a\rangle  \geq c \sum_{k \leq q}   \psi_1(\G_{\a_k}) \Big( | 1- J_K(u)^{\frac{n}{n-2}}\a  _k^{\frac{4}{n-2}}K(a_k) | +  \frac{M_2}{2} (\sum_{r\neq k} \e_{kr} +\frac{1}{\mu_k})\Big) := \ov{\Gamma}_\a .\ee
Such an estimate suggests  to move the variable $\a_i$'s   if $| 1- J_K(u)^{{n}/{n-2}}\a  _i^{{4}/{n-2}}K(a_i) |$ is very large with respect to $\sum_{r\neq k} \e_{kr} +1/\mu_k$.   Furthermore making use of Propositions \ref{p:35} and \ref{pp65}, we derive that
$$ \langle - \n J_K(u), W_a^{in} \rangle \geq c \sum _{i > q} \psi_1(\G_{a_i}) \Big( \frac{| \n K(a_i)|}{\l_i} + \frac{M_2}{2}\Big(\sum_{k\neq i} \e_{ki} +\frac{1}{(\l_i d_i)^{n-2}}+  \frac{1}{\l_i^2}\Big)\Big) := \ov{\G}_a^{in}.$$
Nest  we define $$ W_1:= W_{\L_{in}} + W_{a}^{in} + W_{\a}  + (1/M_2) W_{\L_b}.$$
Using the previous estimates, we obtain
\begin{align}
 \langle - \n J_K(u), W_{1} \rangle & \geq \ov{ \G}_{\L_{in}} + \ov{ \G}_{a}^{in} + \ov{ \G}_{\a}  + (1/M_2) \ov{ \G}_{\L_b} \nonumber \\
& \geq c    \sum_{ i > q }(\psi_1(\G_{\l_i})+\psi_1(\G_{H_i}) + \psi_1(\G_{a_i})) \Big( \sum_{j\neq i} \e_{ij}  +  \frac{H(a_i,a_i)}{\l_i^{n-2}} + \frac{M_2}{2}\frac{1}{\l_i^2} \Big)\nonumber  \\
&  + c \sum _{i > q} \psi_1(\G_{a_i})  \frac{| \n K(a_i)|}{\l_i} + \ov{ \G}_{\a}  + (1/M_2) \ov{ \G}_{\L_b} + O(R_1) .\label{nn1}
\end{align}
Regarding the above estimate, we point that we need to take care of the interaction term $O(\e_{ki})$ contained in the expression    $\ov{ \G}_{\L_b}$.
To that aim, we observe that,  if $\G_{H_k} + \G_{a_k} + \G_{\l_k} \geq 6 $, then the $\e_{ki}$ appears in the lower bound in \eqref{nn1} and therefore we are able to remove the $(1/M_2) \e_{ki}$ by taking  $M_2$ large. But, if  $\G_{H_k} + \G_{a_k} + \G_{\l_k} \leq 6 $, it follows that (see the second assertion of Remark \ref{ree1}) $a_k$ is close to a critical point  $y$ of $K$ and therefore we get $\e_{ki} = O( 1/\l_k^{n-2} + 1/\l_i ^{n-2})$ which is small with respect to our lower bound.\\
Since we are in $V_1(M_0)$, there exists  at least one index $i> q$ such that $\psi_1(\G_{\l_i})+\psi_1(\G_{H_i}) + \psi_1(\G_{a_i}) \geq 1$. This implies that $1/\l_i^2= 1/\mu_i$ appears in the lower bound of \eqref{nn1}. Since all the $\mu_j$'s are of the same order, we are able to make appear all the $1/\mu_j$'s in this lower bound and Lemma \ref{chv1} follows.
\end{pf}

In the next lemma we construct a pseudogradient in the set $V_2(M_0)$. Namely we prove:
\begin{lem}\label{chv2}
There exists a  bounded pseudogradient $W_2$  such  that the following holds : There is a constant $c >  0$ independent of
$u=\sum_{i=1}^{q}\a_i\d_i  + \sum_{i=q+1}^{p+q}\a_i\varphi_i \in V_2(M_0)$ such that the statement of Lemma \ref{chv1}
 holds true with $W_2$ instead of $W_1$.
\end{lem}
\begin{pf}   First, recall that (see Remark \ref{ree1}), in $V_2(M_0)$, each interior concentration point $a_k$  is close to a critical point  of $K$ in $\mathbb{S}^n_+$  and that two interior concentration points $a_i$ and $a_k$ cannot be close to the same critical point which implies that $d(a_i, a_k ) \geq c > 0$  and $\e_{ik} = O( 1/(\l_k\l_i)^{(n-2)/2})$.\\
  Recalling that  $\mathcal{K}^b:= \{ z \in \partial \mathbb{S}^n_+: \n K_1(z)=0\}$ we define the following pseudogradient:
$$ W_2:= W_{\a}  +  W_{\L_b} + \sum _{i \in D_1}    \frac{1}{\l_i}  \frac{\partial \d_i}{\partial a_i} \frac{\n K_1(a_i)}{| \n K_1(a_i)|}  \qquad \mbox{ where } D_1:= \{ i \leq q: d(a_i, \mathcal{K}^b) \geq \eta \}.$$
Using \eqref{cq222}, \eqref{estchalpha} and Proposition \ref{p:34}, we get
\be \langle - \n J_K(u), W_{2} \rangle  \geq  \ov{ \G}_{\a}  +  \ov{ \G}_{\L_b} + \sum _{i \in D_1} \frac{c}{\l_i} + O\Big( \sum_{k \leq q }\big( \frac{1}{\l_i} | \frac{\partial \e_{ki}}{\partial a_i} | +  \l_k d(a_k , a_i ) \e_{ki} ^{\frac{n+1}{n-2}} \big) + R_1^b + \sum_{k > q } \e_{ki} \Big).\label{nn2} \ee
First, taking $i \in D_1$, for $k \leq q$, two cases may occur: $(i)$ either $d(a_i, a_k ) \leq \eta/2$, and in this case we get that $ |\n K(a_k) | \geq c$ and therefore $\mu_k$ and $\l_k$ are of the same order. Thus $\l_i$ and $\l_k$ are of the same order. $(ii)$ or  $d(a_i , a_k ) \geq \eta/2$. In the two cases, we deduce that
$$\e_{ki}= \frac{1+o(1)}{(\l_i \l_k d(a_i , a_k )^2)^{\frac{n-2}{2}}}\, ; \quad  \l_k d(a_k , a_i ) \e_{ki} ^{\frac{n+1}{n-2}}  \leq  \frac{c\, \e_{ki}^\frac{n-1}{n-2} }{\l_i |d(a_i , a_k )} \, \,  \mbox{ and } \, \,  \frac{1}{\l_i} | \frac{\partial \e_{ki}}{\partial a_i} | \leq  \frac{c\,  \e_{ki}}{\l_i d( a_i , a_k )} = o( \e_{ki}) .$$
Secondly, for $i \in D_1$, we have $| \n K_1(a_i) | \geq c(\eta)$ and therefore $\l_i$ and $\mu_i$ are of the same order. Since all the $\mu_j$'s are assumed to be of the same order, we are able to make appear all the $1/\mu_j$'s in the lower bound of \eqref{nn2}.
Finally, for $j\notin D_1$, $(i)$ either
$\G_{\a_j} \geq 2$, in this case, the $|1-J_K(u)^{{n}/(n-2)}\a_j^{{4}/(n-2)} K(a_j)| + \sum \e_{kj}$  appears in $\ov{\G}_{\a}$, $(ii)$ or $\G_{\a_j} \leq 2$ and $\G_{\l_j} \geq 2$, in this case $\sum_{kj} \e_{kj}$ appears in $\ov{\G}_{\L_b}$, $(iii)$ or $\G_{\a_j} + \G_{\l_j} \leq 4$, in this case we are able to make appear   $|1-J_K(u)^{{n}/(n-2)}\a_j^{{4}/(n-2)} K(a_j)| + \sum \e_{kj}$ from $1/\mu_j$. Hence the lemma follows.
\end{pf}

Next  we consider   the third set $V_3(M_0)$. We  notice  that in this subset  each concentration point $a_i$ is close to some critical point of $K$ or $K_1$ and for a critical point  $z \in \partial \mathbb{S}^n_+$ of $K_1$ (resp. $y\in \mathbb{S}^n_+$ of $K$), we denote by
$$B_z :=\{ i \leq q: a_i \mbox{ is close to } z\} \quad ; \quad B_y :=\{ i > q: a_i \mbox{ is close to } y\}.$$
We observe that it follows  from Remark \ref{ree1}  that $\# B_y \leq 1$ for each critical point  $y$ in $\mathbb{S}^n_+$. However, it is possible to have $\#B_z \geq 2$ for some critical points  $z$'s in $\partial \mathbb{S}^n_+$. \\
Next we divide the set $V_3(M_0)$ into  four subsets. The first three ones are defined as follows:
\begin{align*}
   V_3^1:= & \{ u\in V_3(M_0): \exists \, z \mbox{ with } \partial K/\partial \nu(z) =0 \mbox{ and } \# B_z \geq 2\}  , \\
 V_3^2:= & \Big( \{ u: \exists \, z \mbox{ with } \partial K/\partial \nu(z) < 0 \mbox{ and }  B_z \neq \emptyset\} \bigcup \{ u: \exists \, y \mbox{ with } \D K > 0 \mbox{ and } B_y \neq \emptyset\}\\
 & \bigcup \{ u : \exists \, z \mbox{ with } \partial K/\partial \nu(z) = 0 \, \, ; \, \, \D K(z) >  0 \mbox{ and } \# B_z \neq 0\} \Big) \bigcap (V_3(M_0) \setminus V_3^1), \\
  V_3^3:= & \{ u\in V_3(M_0): \exists \, z \mbox{ with } \partial K/\partial \nu(z) > 0 \mbox{ and } \# B_z \geq 2\} \bigcap (V_3(M_0) \setminus  (V_3^1 \cup V_3^2))
\end{align*}
where $y$  is an interior critical point of $K$ and $z$  is a critical point of  $K_1$, and the last one is defined as:

\begin{align}\label{setW}
 \mathcal{W}:=  &  \{   u \in V_3(M_0): \forall i \leq q, a_i \mbox{ is close to   } z_i \in \partial \mathbb{S}^n_+, \mbox{with } \# B_{z_{i}} =1;   ( \partial K/\partial \nu =0 \,  \& \,  \D K < 0)  \nonumber \\
  &  \mbox{  or } \partial K/\partial \nu > 0 \} \bigcap  \{ u \in V_3(M_0), \forall j > q, a_j \mbox{ is close to   } y_j \in \mathbb{S}^n_+, \mbox{with } \# B_{y_{j}} =1 \mbox{ and } \D K(y_j) < 0 \}.
\end{align}

In the next lemma we construct a pseudogradient in the first subset. Namely we prove the following lemma:

\begin{lem}\label{chv31}
There exists a  bounded pseudogradient $W_3^1$  such  that the following holds : There is a constant $c >  0$ independent of
$u=\sum_{i=1}^{q}\a_i\d_i  + \sum_{i=q+1}^{p+q}\a_i\varphi_i \in V_3^1$ such that the statement of Lemma \ref{chv1}
 holds true with $W_3^1$ instead of $W_1$.
\end{lem}
\begin{pf}
Let $z$ be such that $\partial K/\partial \nu(z) = 0$ and $\# B_z \geq 2$. Firstly, we claim that:
\be\label{3333} \mbox{ There exists } k \in B_z  \mbox{ such that: }  \frac{ | \n K_1 (a_k) | }{\l_k } \geq \frac{  M_2 }{\l_k^2} + \frac{1}{ M_2^2} \sum_{j\neq k}  \e_{jk}.\ee
Indeed  arguing by contradiction, we assume that this claim does not hold. Thus, since $z$ is a non-degenerate critical point  of $K_1$, we obtain, for each $k \in B_z$,
 $$c \frac{d( a_k, z )}{\l_k} \leq \frac{| \n K_1(a_k) |}{\l_k} \leq \frac{  M_2 }{\l_k^2} + \frac{1}{ M_2^2} \sum_{j\neq k}  \e_{jk} \leq \frac{  M_2 }{\l_k^2} + \frac{c}{ M_2} \frac{1}{\mu_k} \leq
 \frac{c}{M_2} \frac{| \n K(a_k) |}{\l_k} + c \frac{M_2}{\l_k^2}  \leq c \frac{ d(a_k ,z ) }{M_2\l_k }  + c \frac{M_2}{\l_k^2}$$
 which implies that $\l_k d(a_k , z ) $ is bounded. In addition, from the definition of $\mu_k$, we get
 $$ \frac{1}{\l_k^2} \leq \frac{1}{\mu_k} := \frac{| \n K(a_k) | }{\l_k} + \frac{1}{\l_k^2} \leq c \frac{ \l_k d(a_k , z ) }{\l_k^2} + \frac{1}{\l_k^2} \leq  \frac{c}{\l_k^2}.$$
 Thus, $\mu_k$ and $\l_k^2$ are of the same order for each $k \in B_z$.\\
 Next let $i$ and  $j$ be two different indices in $ B_z$. We deduce that $\l_j$ and $\l_i$ are of the same order and $\l_k d(a_i , a_j )$ is bounded for $k=i,j$. These give a contradiction with the fact that $\e_{ij}$ is small. Hence our claim follows.\\
 \bigskip\noindent
 Furthermore  observe that, for $k$ satisfying \eqref{3333}, it holds that $\l_k d(a_k, z ) \geq c M_2$.\\
 Now, in this region, we define the following vector field:
  $$W_3^1:=  \sum _{i \in D_2}    \frac{1}{\l_i}  \frac{\partial \d_i}{\partial a_i} \frac{\n K_1(a_i)}{| \n K_1(a_i)|}  \qquad \mbox{ where } D_2:= \{ i \leq q: \mbox{ \eqref{3333} holds with }k=i \}.
  $$
Using Proposition \ref{p:34}, we get
\be \langle - \n J_K(u), W_{3}^1 \rangle  \geq c \sum _{i \in D_2} \frac{| \n K_1(a_i) | }{\l_i} + O\Big( \sum_{k \leq q }\big( \frac{1}{\l_i} | \frac{\partial \e_{ki}}{\partial a_i} | +  \l_k d(a_k , a_i ) \e_{ki} ^{\frac{n+1}{n-2}} \big) + R_1^b + \sum_{k > q } \e_{ki} \Big).\label{nn3} \ee
Recall that (see Remark \ref{ree1}), in $V_3(M_0)$, each concentration point $a_k$, for $k > q$ is close to a critical point  of $K$ in $\mathbb{S}^n_+$ which implies that $d(a_i,a_k ) \geq c > 0$ for each $i \leq q$. Hence we get $\e_{ik} = O( 1/(\l_k\l_i)^{(n-2)/2}) $. \\
Moreover for  $i\in D_2$ and  $k \leq q$ with $k\neq i$, two cases may occur:
 (i) either $ \l_k \leq  M_0 ^2M_2^2\l_i$, and in this case we get
$$\frac{1}{\l_i} | \frac{\partial \e_{ik}}{\partial a_i}  | + \l_k d(a_k, a_i ) \e_{ki} ^{\frac{n+1}{n-2}}  \leq c\, \l_k  d(a_i , a_k ) \e_{ik}^{\frac{n}{n-2}}  \leq c\, M_0 M_2\sqrt{\l_k \l_i } d( a_i, a_k ) \e_{ik}^{\frac{n}{n-2}}  \leq  c\, M_0M_2 \e_{ik}^{\frac{n-1}{n-2}}  ,$$
or $(ii)$ $ \l_k \geq  M_0 ^2 M_2^2\l_i$. In this case, since $\mu_k \leq 2 M_0 \mu_i$ and $z$ is a non-degenerate critical point  of $K_1$, it follows that
$$c \frac{d( a_i , z )}{\l_i} \leq \frac{| \n K_1 (a_i) |}{\l_i} \leq \frac{| \n K(a_i) |}{\l_i} + \big(\frac{ 1 }{\l_i^2} - 2\frac{ M_0 }{\l_k^2}\big) \leq 2 M_0 \frac{| \n K(a_k) |}{\l_k} \leq c M_0 \frac{d(a_k , z)}{\l_k}$$
which implies that $ d(a_i, z ) / d(a_k , z ) \leq c M_0 \l_i/\l_k \leq c/ (M_0 M_2^2)$. Thus we deduce that $d( a_i, a_k ) \geq c M_0 M_2^2 d( a_i, z )$. Therefore we obtain
$$\frac{1}{\l_i} | \frac{\partial \e_{ik}}{\partial a_i}  | + \l_k d(a_k, a_i ) \e_{ki} ^{\frac{n+1}{n-2}}  \leq c\, \l_k d(a_i, a_k ) \e_{ik}^{\frac{n}{n-2}}  \leq  \frac{c}{\l_i d(a_i , a_k )} \e_{ik} \leq  \frac{1}{M_0 M_2^2}\frac{c}{\l_i d(a_i, z )} \e_{ik} \leq
\frac{c }{ M_2^3 M_0 } \e_{ik}$$
where we have used the fact that $\l_i d(a_i, z ) \geq c M_2$. Thus \eqref{nn3} becomes
\be\label{555} \langle - \n J_K(u), W_{3}^1 \rangle  \geq c \, \sum _{i \in D_2}  \frac{ | \n K_1(a_i) | }{\l_i} + \frac{ M_2 }{\l_i^2} + \frac{1}{M_2^2} \sum_{j \neq i} \e_{ij} + O(R_1^b + \sum \frac{1}{\l_j ^{n-2}} ).\ee
 Finally, we notice that $ | \n K(a_i) | \leq c d( a_i, z ) \leq c  | \n K_1(a_i) | \leq c  | \n K(a_i) |$. Thus, in \eqref{555}, we can make appear $1/\mu_i$ for $i\in D_2$ and therefore all the $1/\mu_j$'s (since there are of the same order) and the proof follows as the proof of the previous lemmas.
\end{pf}

\begin{lem}\label{chv32}
There exists a  bounded pseudogradient $W_3^2$  such  that the following holds : There is a constant $c >  0$ independent of
$u=\sum_{i=1}^{q}\a_i\d_i  + \sum_{i=q+1}^{p+q}\a_i\varphi_i \in V_3^2$ such that the statement of Lemma \ref{chv1}
 holds true with $W_3^2$ instead of $W_1$.
\end{lem}
\begin{pf}  Let $D_1:= \cup_{y: \D K(y)  > 0} B_y$,  $D_2:= \cup_{z: \partial K/\partial \nu(z) < 0} B_z$  and   $D_3:= \cup_{z: \partial K/\partial \nu(z) = 0 \, ;\,  \D K(z)  > 0} B_z$. We  divide this region into two subsets:

{\it 1st subset:} If $D_1 \cup D_2 \neq \emptyset$. In this case, we define
$$W_3^{21}:= - \sum_{i\in D_1 \cup D_2} \l_i \frac{\partial \varphi_i}{\partial \l_i} .$$
By using the first assertion of Remark \ref{ree1} and Propositions \ref{p:33} and \ref{pp65}, it follows that
$$ \langle - \n J_K(u), W_{3}^{21} \rangle  \geq c \sum _{i \in D_1 \cup D_2} \Big( \sum_{ j\neq i} \e_{ij} + \frac{1}{\mu_i} + O\big( \sum \frac{1}{\l_j ^{n-2} }+ R_1^b + R_1\big) \Big).$$
Hence, the proof follows.

{\it 2nd subset:} $D_3 \neq \emptyset$. Note that, since we are outside of $V_3^1$,  for $i \in B_z$ with $\partial K/\partial \nu(z) = 0$, it holds that $B_z =\{ i \}$, that is $ d( a_i, a_j ) \geq c > 0$ for each $j \neq i$. We define
$$ W_3^{22} := \sum_{i \in D_3} \psi_1(\l_i | \n K_1(a_i) | /M) \frac{1}{\l_i} \frac{\partial \d_i}{\partial a_i} \frac{  \n K_1(a_i)  }{| \n K_1(a_i) |} -  \l_i \frac{\partial \d_i}{\partial \l_i} $$ where $M$ is a large constant. We point out  that $W_3^{22}$ is exactly the sum of  of the vector fields $\underline{W}_1^{3}$ (defined in \eqref{undW13}) with $\gamma=-1$. Furthermore, the presence of the function $\psi_1$ implies that the point $a_i$ moves only if $ | \n K_1(a_i) | \geq M/\l_i$.\\
 Using Propositions \ref{p:33} and \ref{p:34}, we get
\begin{align}  \langle - \n J_K(u), W_{3}^{22} \rangle  \geq & c \sum _{i \in D_3}    \psi_1(\l_i | \n K_1(a_i) | /M) \Big( \frac{ | \n K_1(a_i) | }{\l_i} + \frac{1}{\l_i ^2} \Big) \label{zzaa1}\\
&  -  \Big(  \frac{c_3}{\l_i} \frac{\partial K}{\partial \nu} (a_i) - c \frac{\D K(a_i)}{\l_i ^2}\Big) + O\Big( \sum \frac{1}{\l_j ^{n-2} } +R_1^b \Big) \nonumber \end{align}
which has the same form as  \eqref{mm1}. Hence, the same computations and arguments hold and the proof of the lemma follows.
\end{pf}

\subsubsection{Ruling out collapsing phenomena  }\label{s:simpleblow1}

We point that, the main difference between the $\mathbb{S}^n$-case (or the case of an interior blow up point for the $\mathbb{S}^n_+$-case) and the boundary blow up point case  relies essentially on the behavior of  the leading term   in Propositions \ref{p:33} and \ref{pp65} (namely  the $\l$-term). Indeed when  $\partial K/\partial\nu (z) \ne 0$ and $a_i$ is close to  a boundary critical point $z\in \partial \mathbb{S}^n_+$, the leading term behaves  like $c/\l_i$, while  for the $\mathbb{S}^n$-case (or the case of an interior blow up point in  the $\mathbb{S}^n_+$-case), for  $a_i$ close to an interior critical point $y$  with $\D K(y) \ne 0$, this term behaves  like $c/\l_i^2$. This difference on the behavior of the leading term plays a crucial role in the nature of the \emph{critical point at infinity}. Indeed in \cite{AB}, for $z$ a critical point  of $K_1$ (which is not  local maximum) satisfying  $\partial K/\partial\nu (z) > 0$, we proved that $z$ is not a simple blow up point in the sense that $B_z $ contains more than one concentration point. In the following lemma, we  consider the case of a  local maximum point of $K_1$ satisfying  $\partial K/\partial\nu (z) > 0$ and we will prove  that $z$ is a simple blow up point.
Namely we prove
\begin{lem}\label{simple}
Let $z$ be a  non degenerate local maximum   of $K_1$  with $\partial K/\partial\nu (z) > 0$. Then $z$ is a \emph{simple blow up}. More precisely if $\# B_z := \, \# \{ a_i; \mbox{close to }  z\} \, :=  q_1 \geq 2$, then
 $J_K$ admits in the set $V(q_1,q_1,0,\e)$ \emph{a compactifying  bounded pseudogradient} $W(z,q_1)$. Namely there exits  a constant $c>0$ independent of
$u=\sum_{i=1}^{q_1}\a_i\d_i $ such that
$$ \langle-\nabla J_K(u),W(z,q_1)\rangle \geq  c \,  \sum_{i\leq q_1}  \Big(\frac{1}{\l_i^{2-{1}/{(n-2)}}}+|1-J_K(u)^{\frac{n}{n-2}}\a_i^{\frac{4}{n-2}} K(a_i)|^{2-\frac{1}{n-2}}  \Big) + c\, \sum_{k\neq r} \e_{kr}^{\frac{n-1}{n-2}}. $$
Furthermore, the concentration rates $\l_i$'s  do not increase  along the flow lines generated by this pseudogradient.
\end{lem}

For the proof of Lemma \ref{simple}, we make use of the following technical results.

\begin{lem}\label{dereps}
Let $a_i, a_j \in \partial \mathbb{S}^n_+$ be concentration points  such that the corresponding rates  $\l_i$ and $\l_j$ are of the same order and $d(a_k, b) \to 0$ for $k=i,j$ for some point $b\in \partial \mathbb{S}^n_+$. Then we have
$$ e_{ij}:=  \frac{ \partial \e_{ij} }{\partial a_i } (b- \langle a_i , b\rangle a_i) + \frac{ \partial \e_{ij} }{\partial a_j } (b- \langle a_j , b\rangle a_j)  \geq  c \, \e_{ij}. $$
\end{lem}

\begin{pf}
Easy computation implies that $$ \frac{ \partial \e_{ij} }{\partial a_i } = (n-2) \l_i \l_j (a_j - a_i ) \e_{ij}^{n/(n-2)}.$$ Thus we get
\begin{align*}
e_{ij} & =  (n-2) \l_i \l_j  \e_{ij}^{n/(n-2)}\big( \langle a_j - a_i , b - \langle a_i , b\rangle a_i\rangle + \langle a_i - a_j , b - \langle a_j , b\rangle a_j\rangle \big)\\
& =  (n-2) \l_i \l_j  \e_{ij}^{n/(n-2)} \langle a_j + a_i , b \rangle (1- \langle a_i , a_j \rangle )\\
& =  (n-2) \l_i \l_j  \e_{ij}^{n/(n-2)}| a_i - a_j |^2 (1+o(1))  \, \geq \,  c \e_{ij}.
\end{align*}
 where $| a_i - a_j | $ is the euclidian norm of $a_i-a_j$ seen as a vector in $\R^{n+1}$.
\end{pf}


\begin{lem}\label{derK}
Let $a, \, h \in \partial \mathbb{S}^n_+$ be close to a non degenerate local maximum  $z$ of $K_1$. Then it  holds that
$$ \frac{1}{K_1(a)^{n/2}} \n K_1(a) \big(h - \langle a, h \rangle a\big) \geq  - \frac{1}{K_1(h)^{n/2}} \n K_1(h) \big(a - \langle a, h \rangle h\big)  + c\,  | a - h |^2.$$
\end{lem}
\begin{pf}
Let $$ \beta(t):= \frac{ h + t(a-h)}{| h + t(a-h) |} \quad , \quad g(t):= \frac{ 2/(n-2) }{ K_1(\beta(t))^{(n-2)/2}} \quad \mbox{ for } t\in [0,1].$$
It is easy to get that
$$ \beta'(t) = \frac{1} { | h + t(a-h) |} \Big( a - h - \langle \beta(t) , a-h \rangle \beta(t) \Big) \quad , \quad  \langle \beta(t) , a-h \rangle = O(| a - h |^2),$$
and therefore it holds that $| \beta'(t) | = |a-h | (1+o(1))$ uniformly in $t\in [0,1]$. Furthermore, easy computations imply that $| \beta ''(t) | = O( |a-h | ^2)$ uniformly in $t\in [0,1]$. In another hand, we have
$$ g'(t) = \frac{-1}{  K_1(\beta(t))^{n/2}} \n K_1(\beta(t)) \big(\beta'(t)\big) $$
and, since $a$ and $h$ are close to a non degenerate maximum critical point  $z$ of $K_1$, we derive that
$$g''(t) = o( | \beta'(t) |^2 ) - \frac{1}{  K_1(\beta(t))^{n/2}} D^2 K_1(\beta(t)) \big(\beta'(t), \beta'(t)\big) + o( | \beta''(t) | ) \geq c | a- h |^2 \,\,  (\mbox{uniformly in }  t\in [0,1]).$$
Now,
$$ \frac{1}{K_1(a)^{n/2}} \n K_1(a) \big(h - \langle a, h \rangle a\big) + \frac{1}{K_1(h)^{n/2}} \n K_1(h) \big(a - \langle a, h \rangle h\big) = g'(1) - g'(0) = \int _0^1 g''(t) \, dt $$ which implies the lemma.
\end{pf}


\begin{pfn}{\bf of Lemma \ref{simple}}
For the construction of a suitable vector field satisfying the properties required in Lemma \ref{simple} as well for later purposes we will use some constants $M_0$, $M_2$ and $M_4$ which are required to be large and to satisfy
\be\label{M0M2M4} \frac{ M_0}{M_4 ^2} \mbox{ small } \quad , \quad \max \Big( \frac{ M_2}{M_0^{1/(q+p-1)}}  \, ;\,  \frac{ M_2^{(n-1)/(n-2)}}{M_0^{(1/2 + 1/(n-2))/(q+p-1)}} \Big) \mbox{ small }.\ee
The first requirement  is used in \eqref{ee1} and \eqref{eee1} below while the second one is used when studying a remainder term of \eqref{M2M0} and the last one is used in \eqref{M0M2} in the proof of Lemma \ref{chv4}. \\
In view of the pseudogradient constructed in Lemmas \ref{chv2} and \ref{chv4}, it is enough to construct a pseudogradient satisfying the above estimate in the following set:
  $$ V(z,q_1,\eta, \e,M_0):= \{ u\in V(q_1,q_1,0, \e): \l_{\max} \leq M_0 \l_{\min} \, ; \,  d( a_i, z ) < \eta \, \, ; \G_{\l_i} \leq 2 \, \mbox{ and } \G_{\a_i} \leq 2 \, \, \forall \, i\}.$$
 Moreover,  since the $\l_i$'s are of the same order, we have that $\e_{ij}=(1+o(1))/(\l_i \l_j d(a_i,a_j)^2)^{(n-2)/2}$ and therefore $d( a_i ,a_j ) \geq c/\l_1^{(n-3)/(n-2)}$  for each $i\neq j$ (since $\G_{\l_i}$ is bounded). We want to construct a pseudogradient which moves the concentration points $a_i$ to their barycenter  and prove that along its flow lines the Euler-Lagrange functional decreases. To this aim, let $i$ and $i_1$ be such that $d( a_i , a_{i_1} ) := \min d(a_r, a_\ell )$ and define $L_i:= \{i,i_1\}$. Next  let $M_4$ be a large positive constant,  for such an index $i$, we define inductively a sequence $L_{i}^s$ by setting
\begin{align*} & L_i^1 := \{ j : \, \exists \, \, \ell \in L_i \, \, s.t. \, \, d(a_j , a_\ell ) \leq M_4  d( a_i , a_{i_1} ) \}  \quad \mbox{ and} \\
&  L_i^s := \{ j : \, \exists \, \, \ell\in L_i^{s-1} \, \, s.t. \, \, d(a_j, a_\ell ) \leq M_4 \max_{r, t \in L_i^{s-1} } d( a_r , a_t ) \} .\end{align*}
Observe that, since we have only $q_1$ points and $\# L_i = 2$, then there exists $m \leq q_1-1$ such that $L_i^{m+1}= L_i^m$ and we set
$L_i^* := L_i^m$ where $m$ is the first index  such that $L_i^{m+1}= L_i^m$. We remark that $L_i \subset L_i^*$.
Next we want to move the points $a_j$'s, for $j \in L_i^*$, to their center of mass. For this aim, let  ${\bf \ov a_i}$ be defined as
\be\label{barai1}  {\bf \ov a_i} := \frac{b_i}{| b_i| } \quad \mbox{ where } \quad b_i \in \R^{n+1} \mbox{ satisfying }  \sum_{j\in L_i^*} (b_i - a_j) = 0.\ee
Note that, it is easy to see that ${\bf \ov a_i}$ satisfies
\be\label{barai2}  {\bf \ov a_i} \in \partial \mathbb{S}^n_+ \quad \mbox{ and } \quad  \sum_{j\in L_i^*} a_j - \langle a_j , {\bf \ov a_i} \rangle {\bf \ov a_i} = 0.\ee
Now we define  the following vector field:

$$ W_3^i := \frac{1}{\l_i \gamma_i } \sum _{j\in L_i^*} \a_j\frac{\partial \d_j}{\partial a_j} ( {\bf \ov a_i} - \langle a_j , {\bf \ov a_i} \rangle   a_j ) \qquad \mbox{ where } \quad \gamma_i := \max_{j\in L_i^*} d(a_i, a_j).$$
We note that $L_i^*$ has two important  properties:
\begin{itemize}
\item If $k,\ell  \in L_i^*$, we have $d( a_k ,a_\ell ) \leq c M_4^m  d( a_{i} , a_{i_1} )$.
\item  If $k \notin L_i^*$, then, for each $j \in L_i^*$,  we have  $d( a_j, a_k ) \geq M_4 \max_{r,\ell \in L_i^*} d(a_r, a_\ell ) $. Hence, for $k \notin L_i^*$ and $ j \in L_i^*$, choosing ${M_0^{(n-2)/2}}/{M_4^{n-2}}$ small, it follows that for every  $ \ell \in L_i^*$, we have that:
\end{itemize}
\begin{align}  & | \frac{\partial \e_{jk}}{\partial a_j} | | {\bf \ov a_i} - \langle a_j , {\bf \ov a_i} \rangle   a_j |   \leq \frac{c d({\bf \ov a_i}, a_j)  }{(\l_j \l_k) ^{\frac{n-2}{2}} d( a_j, a_k )^{n-1}}   \leq \frac{ M_0^{(n-2)/2}}{M_4^{n-1}}\frac{c} {(\l_j \l_\ell) ^{\frac{n-2}{2}} d(a_j,a_\ell )^{n-2}} = o\Big( \e_{j\ell } \Big)\label{ee1}\\
& \e_{jk} \leq  \frac{c }{(\l_j \l_k) ^{(n-2)/2} d(a_j, a_k )^{n-2}} \leq  \frac{c M_0^{(n-2)/2}}{M_4^{n-2}}\frac{1} {(\l_j \l_\ell) ^{(n-2)/2} d( a_j, a_\ell )^{n-2}} = o\Big( \e_{j\ell } \Big) \label{eee1}\end{align} (by using  \eqref{M0M2M4}).
We note that, in this region, we have $|1- J_K(u)^{\frac{n}{n-2}}\a  _j^{\frac{4}{n-2}}K(a_j)| \leq c M_2/{\l _j}$ for each $j$, hence Proposition  \ref{p:34} can be written as :
\begin{align}
\langle & \nabla J_K(u),   \a_j\frac{\partial \d _j}{\partial a_j}\rangle  = \l_j  \left[ c_4\left(1- J_K(u)^{\frac{n}{n-2}}\a _i ^{\frac{4}{n-2}}K(a_i)\right)+ J_K(u)^{\frac{n}{n-2}}\a _i ^{\frac{4}{n-2}}\frac{c_5}{\l_i}\frac{\partial K}{\partial\nu}(a_i)\right] e_n  \label{devda2}\\
& - J_K(u)c_2 \sum_{k\neq j} \a _j\a_k  \frac{\partial\e _{kj}}{\partial a_j}   - 8 c_5 \, J_K(u)^{-\frac{n-2}{2}} \frac{1}{K(a_j)^{n/2}} \n K_1(a_j) + O\Big(\frac{1}{\l} + \l \sum \e _{kr }^{\frac{n}{n-2}}\ln (\e _{kr}^{-1})\Big).\nonumber
\end{align}
Hence we derive  that:
\begin{align}
\langle -\nabla J_K(u), & W_3 ^i \rangle =   \frac{J_K(u)c_2}{\l_i \g_i} \sum_{k\neq j; j\in L_i^* } \a _j\a_k  \frac{\partial\e _{kj}}{\partial a_j}( {\bf \ov a_i} - \langle a_j , {\bf \ov a_i} \rangle   a_j )\nonumber\\
&   + \frac{8 c_5 \, J_K(u)^{(2-n)/2} }{\l_i \g_i} \sum_{j\in L_i^*}\frac{1}{K(a_j)^{n/2}} \n K_1(a_j) ( {\bf \ov a_i} - \langle a_j , {\bf \ov a_i} \rangle   a_j )  + O\Big( \frac{1}{\l^2}   + \sum \e _{kr }^{\frac{n}{n-2}}\ln (\e _{kr}^{-1})\Big). \label{ee2}
\end{align}
Next we notice  that, using Lemma \ref{dereps}, il holds
\be\label{ee3}\frac{\partial\e _{kj}}{\partial a_j}( {\bf \ov a_i} - \langle a_j , {\bf \ov a_i} \rangle   a_j )  + \frac{\partial\e _{kj}}{\partial a_k} ( {\bf \ov a_i} - \langle a_k , {\bf \ov a_i} \rangle   a_k )\, \geq c \, \e_{kj}, \qquad \mbox{ for each } k, j \in L_i^*.\ee
Furthermore, using Lemma \ref{derK} (with $h={\bf \ov{a}_i}$), it holds that
\begin{align*}
 \sum_{j\in L_i^*}\frac{1}{K(a_j)^{n/2}} \n_T K(a_j) ( {\bf \ov a_i} - \langle a_j , {\bf \ov a_i} \rangle   a_j ) &  \geq  \sum_{j\in L_i^*}\frac{-1}{K( {\bf \ov{a}_i} )^{n/2}} \n_T K( {\bf \ov{a}_i} ) ( a_j  - \langle a_j , {\bf \ov a_i} \rangle   {\bf \ov{a}_i} )  + c  \sum_{j\in L_i^*} | a_j - {\bf \ov{a}_i}| ^2 \\
 & \geq c  \sum_{j\in L_i^*} | a_j - {\bf \ov{a}_i}| ^2 \quad (\mbox{since }  \sum_{j\in L_i^*}  a_j  - \langle a_j , {\bf \ov a_i} \rangle   {\bf \ov{a}_i} = 0\,  (\mbox{see \eqref{barai2}})).
\end{align*}
 Thus we get
\be\label{aaa1}
\langle -\nabla J_K(u), W_3^i \rangle \geq  c \sum_{k, j \in L_i^*} \frac{\e_{kj}}{\l_i \gamma_i} + \sum_{j\in L_i^*} \frac{ d( a_j,{\bf \ov{a}_i})^2 }{\l_i \gamma_i} +O\Big( \sum \e_{\ell r}^{\frac{n}{n-2}} \ln \e_{\ell r} ^{-1} + \frac{1}{\l_i ^2}\Big).
\ee
Now, since $\gamma_i:= \max_{k,r \in L_i^*} d(a_k,a_r)$ is of the same order of all the $d(a_\ell, a_j )$'s, we derive that ${\e_{kj}}/{\l_i \gamma_i} \geq c \e_{kj} ^{(n-1)/(n-2)}$. Furthermore, $\sum _{j\in L_i^*}  d(a_j ,{\bf \ov a_i})^2 \geq c \sum _{j, r \in L_i^*}  d( a_j, a_r )^2 $ and therefore
$$\sum _{j\in L_i^*}  d(a_j, {\bf \ov a_i})^2 / (\l_i \gamma_i) \geq \sum _{j, r \in L_i^*}  d(a_j,  a_r ) / \l_i \geq c /\l_i ^{2- 1/(n-2)}.$$
Hence, in the lower bound of \eqref{aaa1}, we are able to make appear $1 /\l_i ^{2- 1/(n-2)}$ and therefore (since all the $\l_j$'s are of the same order and $\G_{\a_k} \leq 4$ for each $k$) we are able to make appear all the $1 /\l_j ^{2- 1/(n-2)} $'s and $| 1- J_K(u)^{n/(n-2)} \a_j^{4/(n-2)} K(a_j) |^{2- 1/(n-2)} $'s. Concerning the  $\e_{kr}$, we note that the $\e_{kj}$'s which appeared in the lower bound, are only for the indices $k, j \in L_i^*$. Hence we need to make appear $\e_{jr}$ for $j\notin L_i^*$. For this aim, we remark that, for each $j, \ell$, we have $ d( a_j, a_\ell ) \geq  d( a_i, a_{i_1} )$ (by the definition of $i$ and $i_1$), in addition we have that the $\l_k$'s are of the same order. Hence we deduce that $\e_{ii_1} \geq c \e_{j \ell} $. Hence the proof of the lemma follows.
\end{pfn}




In the next lemma we rule out non simple blow up for a \emph{ mixed configuration } involving local  maxima on the boundary and other interior blow up points. Namely we prove:


\begin{lem}\label{chv33}
There exists a  bounded pseudogradient $W_3^3$  such that the following holds : There is a constant $c >  0$ independent of
$u=\sum_{i=1}^{q}\a_i\d_i  + \sum_{i=q+1}^{p}\a_i\varphi_i \in V_3^3$ such that
$$ \langle-\nabla J_K(u),W_3^3 \rangle \geq     \sum_{i=1}^{q+p} \frac{c}{\mu_i^{\frac{2n-5}{n-2}}}+ c \sum_{i=1}^{q}|1-J_K(u)^{\frac{n}{n-2}}\a_i^{\frac{4}{n-2}} K(a_i)|^{\frac{2n-5}{n-2}}  + c\, \sum_{k\neq r} \e_{kr} ^\frac{n-1}{n-2} + c \sum_{i=q+1}^{q+p}\big( \frac{| \n K(a_i)|}{\l_i}\big)^{\frac{2n-5}{n-2}}  $$
Furthermore, the $\l_i$'s do not increase  along the flow lines generated by the pseudogradient $W^3_3$.
\end{lem}

\begin{pf} Let $z_1, \cdots ,z_\ell$ be the critical points  of $K_1$ satisfying  $\partial K /\partial \nu (z_j) > 0$ and $\# B_{z_j} \geq 2$. We decompose  $u$  as follows:
$$u:=\sum_{i=1}^\ell u_i + u_{\ell +1} \quad \mbox{ where } u_i:= \sum_{k\in B_{z_i}} \a_k \d_k \mbox{ and } u_{\ell +1} := u-  \sum_{i=1}^\ell u_i .$$
From the second and the third assertions of Remark \ref{ree1}, it follows that each concentration point $a_j$ of $u_{\ell +1}$ satisfies $ | a_j - a_k | \geq c $ for each $k \neq j$ and it is close to a critical point  of $K_1$ with $\partial K/\partial \nu \geq 0$ or a critical point  of $K$  in $\mathbb{S}^n_+$ with $\D K < 0$. Furthermore, for $j \in B_{z_i}$, we have $  | a_j - a_k | \geq c $ for each $k \notin B_{z_i}$. Hence the mutual interaction between two clusters $B_{z_i}$ and $B_{z_j}$  for $i \neq j $ is   negligible with respect to the other terms. In this situation, we define the following vector field
 $$ W_3^3 := \sum_{i=1}^\ell W(z_i, \# B_{z_i})$$
 where  $W(z_i, \# B_{z_i})$ is defined in Lemma \ref{simple}. Hence we obtain
\be \label{666}  \langle-\nabla J_K(u),W_3^3 \rangle = \sum_{i=1}^\ell  \langle-\nabla J_K(u), W(z_i, \# B_{z_i}) \rangle = \sum_{i=1}^\ell   \langle-\nabla J_K(u_i), W(z_i, \# B_{z_i}) \rangle + \sum_{k\in B_{z_i}; j\notin B_{z_i}} O\big( \e_{kj}\big) .\ee
We observe  that, for $k\in B_{z_i}$, we have $\mu_k$ and $\l_k$ are of the same order. Moreover  we are in the case where all the $\mu_j$'s are of the same order. Thus, using Lemma \ref{simple}, we are able to make appear all the $1/\mu_j^{2-1/(n-2)}$'s in the lower bound of \eqref{666} (and therefore all the $|1-J_K(u)^{{n}/(n-2)}\a_i^{{4}/(n-2)} K(a_i)|^{2-{1}/(n-2)}$'s and the $ ({| \n K(a_i)|}/{\l_i})^{2-{1}/(n-2)}$'s (since the $\G_{\a_k}$'s and the $\G_{a_i}$'s are bounded). In addition, for $j\notin B_{z_i}$ and $k\in B_{z_i}$, we have $$\e_{kj} \leq  \frac{c }{(\l_j \l_k)^{(n-2)/2}}  \leq  \begin{cases}
& o(1/\l_k^2) \mbox{ if } n \geq 6 , \\
& c/\l_k^2 + c/\l_j^4  \mbox{ if } n =5.\end{cases} $$
Therefore, our lemma follows from Lemma \ref{simple}.
 \end{pf}

\begin{lem}\label{chv34}
There exists a  bounded pseudogradient $\mathcal{V}$  satisfying the following estimate : \\
There is a constant $c >  0$ independent of
$u=\sum_{i=1}^{q}\a_i\d_i  + \sum_{i=q+1}^{p}\a_i\varphi_i \in \mathcal{W}$ such that
\eqref{estchv1} holds true with    $\mathcal{V}$  instead of $W_1$.\\
Furthermore in the subset of $\mathcal{W}$ such that $ \l_i |\n K_1(a_i) |$ is bounded, the  $\l_i$'s are increasing functions along the flow lines generated by the pseudogradient $\mathcal{V}$.
\end{lem}

\begin{pf}
Let  $\psi_1$ be a $C^\infty$ cut of  function defined by $\psi_1 \in [0,1]$, $\psi_1(t) = 1 $ if $t \geq 2$ and $\psi_1(t)=0$ if $t\leq 1$. \\
We define the following vector field:

$$\mathcal{V}:= W_\a +  W_a^{in} + W_a^b + \sum_{i=1}^{p+q} \l_i \frac{\partial \varphi_i}{\partial \l_i} $$ where  $W_a^b := \sum_{i\in I_b} \psi_1 (\l _i | \n K_1(a_i) | / M_2) (1/\l_i) (\partial \d_i/\partial a_i)(\n K_1(a_i) / | \n K_1(a_i) |)$ and $W_a^{in}$ (resp. $W_\a$) is defined in \eqref{aaz1} (resp. \eqref{aaz2}).

Observing that in $\mathcal{W}$ we have  $\e_{ij} = O( 1/\l_i ^{n-2} + 1/\l_j ^{n-2})$ for each $i \neq j$ and  using  Propositions \ref{p:33}, \ref{p:34}, \ref{pp65}  the lemma follows.
\end{pf}




\subsubsection{ Ruling out bubble towers phenomena}\label{s:simpleblow2}
In this subsection we prove any configuration of points of non comparable  concentration rates is not critical at infinity. Indeed one can construct in the neighborhood of such     points a \emph{compactifying pseudogradient}. Namely we prove that:
\begin{lem}\label{chv4}
There exists a  bounded pseudogradient $W_4$  such  that the following holds : There is a constant $c >  0$ independent of
$u=\sum_{i=1}^{q}\a_i\d_i  + \sum_{i=q+1}^{p+q}\a_i\varphi_i \in V_4(M_0)$ such that
$$ \langle-\nabla J_K(u),W_4 \rangle \geq  c \,  \sum_{i=1}^{q+p} \frac{1}{\mu_i^\frac{n-1}{n-2}}+ c \sum_{i=1}^{q}|1-J_K(u)^{\frac{n}{n-2}}\a_i^{\frac{4}{n-2}} K(a_i)|^\frac{n-1}{n-2}  + c\, \sum_{k\neq r} \e_{kr} ^\frac{n-1}{n-2} + \sum_{i=q+1}^{q+p}\big( \frac{| \n K(a_i)|}{\l_i}\big)^\frac{n-1}{n-2} . $$
Furthermore, $\max \mu_i$ deos not increase  along the flow lines generated by this pseudogradient.
\end{lem}
\begin{pf}
For $u= \sum_{i=1}^{q}\a_i\d_{a_i,\l_i}  + \sum_{i=q+1}^{p+q}\a_i\varphi_{a_i,\l_i}$, we denote
$$
\mathcal{I}_{in} := \{ i=1, \cdots, p+q; \, a_i \in \mathbb{S}^n_+ \} \quad  \& \quad \mathcal{I}_{b} := \{ i=1, \cdots, p+q; \, a_i \in \partial \mathbb{S}^n_+ \} .
$$
Next we  reorder the parameters $\mu_i$'s as: $\mu_1 \leq \cdots \leq \mu_{p+q}$ and define the following subset of indices:
$$I:= \{1\} \cup \{ i \geq 2 : \mu_k \leq M_0^{1/(p+q -1)} \mu_{k-1} \mbox{ for each } k\leq i\}.$$
Since we are in $V_4(M_0)$, we have $\mu_{\max} > M_0 \mu_{\min}$, it follows that $p+q \notin I$.  In this region, we write $u$ as
$$ u := u_1 + u_2 \quad \mbox{ where } u_1 := \sum_{i\in I} \a_i \varphi_i \mbox{ and } u_2 := u - u_1.$$
Let $k_0:= \max I$ (then we have $k_0 < p+q$). It follows that $\mu_{k_0} \leq M_0^{(k_0 - 1 )/ (p+ q -1)} \mu_1 := \ov{ M}_0 \mu_1$,  $\mu_{k_0+1} \geq M_0^{1/(p+q -1)} \mu_{k_0}$ and therefore $u_1 \in V_1(\ov{ M}_0) \cup V_2(\ov{ M}_0) \cup V_3(\ov{ M}_0)$.


Furthermore we introduce the following notation
$$
D^4_1:= \{ i \in \mathcal{I}_{in}: \G_{\l_i} + \G_{a_i} + \G_{H_i} \geq 6\} \quad \& \quad D^4_2:= \{ i \in \mathcal{I}_{b}: \G_{\l_i} + \G_{\a_i} \geq 4\}
$$
and set
$$
i_0:=
\begin{cases}
  \min D^4_1, & \mbox{if }  D^4_1 \neq \emptyset \\
  p+q+1 , & \mbox{otherwise}.
\end{cases}
\qquad
j_0:=
\begin{cases}
  \min D^4_2, & \mbox{if }  D^4_2 \neq \emptyset \\
  p+q+1 , & \mbox{otherwise}.
\end{cases}
$$
Next   we define in case $ D^4_1 \cup D^4_2 \neq \emptyset$ the following vector fields:
$$ W_{i_0}:= - \sum_{i\geq i_0; i\in \mathcal{I}_{in}} 2^i \l_i \frac{\partial \varphi_i}{\partial \l_i} \quad \mbox{ and }\quad  W_{j_0}:= - \sum_{j\geq j_0; i\in \mathcal{I}_{b}} 2^j \l_j \frac{\partial \d_i}{\partial \l_i} $$ and as in the proof of Lemma \ref{chv1}, we define
$$ W_4^0:=  W_{i_0} + (1/M_2) W_{j_0} + W_{\a} + W_a^{in}$$
where  $W_a^{in}$ (resp. $W_{\a}$) is defined in \eqref{aaz1} (resp. \eqref{aaz2}). Following the proof of Lemma \ref{chv1} and using Lemma \ref{epsij1}, we get
\begin{align}
 \langle - \n J_K(u), W^{0}_4 \rangle  \geq  & \, \, \ov{ \G}_{a}^{in}+ c \sum_{i\geq i_0; i\in \mathcal{I}_{in}} \Big( \sum_{\ell \neq i} \e_{i\ell} + \frac{1}{(\l_i d_i)^{n-2}} +O(R_1)\Big)  +  \frac{c}{\mu_{i_0}}\nonumber \\
 &  + \ov{ \G}_{\a}  + \frac{c}{M_2} \sum_{ j\geq j_0; j\in \mathcal{I}_{b}} \Big( \sum_{\ell \neq j } \e_{j\ell }  +O\big(R_1^b +   + c\frac{c}{\mu_{j_0}}+ \sum_{\ell \in \mathcal{I}_{in}} \e_{j\ell }\big)\Big) := \ov{\G}_4. \label{Gamma4}
\end{align}
Observe that, concerning the last term, for $\ell \in \mathcal{I}_{in}$, $(i)$ either $\ell \geq i_0$, then the $\e_{j\ell}$ exists in the second term of  this formula and one takes  $M_2$  large to absorb the last term, or $(ii)$  $\ell < i_0$ and in this case  by Lemma \ref{pointint}, the concentration point $a_\ell$ is close to a critical point  $y$ of $K$ in $\mathbb{S}^n_+$ and then $\e_{j\ell} \leq c(1/ \l_j^{n-2} + 1/ \l_\ell^{n-2})$. Hence, we can in either case absorb the last term. \\
Furthermore we notice that if   $ D^4_1 \cup D^4_2 \neq \emptyset$ and if $i_0 \in I $ or if $j_0 \in I$ then we can include all the indices in $I$ in the lower bound of \eqref{Gamma4}. Otherwise  to make appear the terms corresponding to  these indices  we argue as follows:

{\it Case 1:}  If $u_1  \in V_1(\ov{ M}_0) \cup V_2(\ov{ M}_0) \cup (V_3(\ov{ M}_0)  \setminus V_3^3) $. In this region, we define the following vector field:
$$ W_4^1:= W_4^0 + (1/M_2^2) \wtilde{W}(u_1),$$
where $\wtilde{W}$ is the convex combination  of the pseudogradients constructed in $V_1(\ov{ M}_0)$,   $V_2(\ov{ M}_0)$ and  $V_3(\ov{ M}_0)  \setminus V_3^3$. It follows then that
\begin{align}  \langle - \n J_K(u), W_4^1 \rangle \geq &\, \,  \ov{\G}_4 +  \frac{1}{M_2^2}\Big( \,  \sum_{i \in I} \frac{c}{\mu_i}+ c\sum_{i\in I\cap \mathcal{I}_b}|1-J_K(u)^{\frac{n}{n-2}}\a_i^{\frac{4}{n-2}} K(a_i)| \nonumber\\
 & + c \sum_{k\neq r; k,r\in I} \e_{kr}  + c \sum_{i\in I\cap \mathcal{I}_{in}} \frac{| \n K(a_i)|}{\l_i} + O( \sum _{j \in I; \ell \notin I} \e_{j\ell })\Big). \label{M2M0}\end{align}
To complete the proof, it remains to absorb the last term. To  this aim, we notice that: \\
$(i)$ if "$\ell \in \mathcal{I}_{in}$ with $\ell \geq i_0$ or $\ell \in \mathcal{I}_b$ with $\ell \geq j_0$", then the term $\e_{j\ell}$ is already in $\ov{\G}_4$   the lower bound of \eqref{Gamma4}.  Taking  $M_2$  large, we will be able to absorb this term.\\
 $(ii)$ if  "$\ell \in \mathcal{I}_{in}$ with $\ell < i_0$ or $\ell \in \mathcal{I}_b$ with $\ell < j_0$", then there holds:  $\e_{j\ell} \leq c \frac{M_2}{\mu_\ell} \leq c (M_2/ { M}_0^{1/(q+p-1)})\frac{1}{\mu_{k_0}} = o(1/\mu_{k_0})$ by choosing $M_2/ { M}_0^{1/(q+p-1)}$ small enough (see \eqref{M0M2M4}) and where $k_0 := \max I$. Hence, we are also able to remove this term. (Recall that, in Lemmas \ref{chv1}-\ref{chv32}, \ref{chv34}, the constant over $\mu_{max}$ is independent of $M_0$ and $M_2$).  Hence the estimate in the first case follows as in the proof of the previous lemmas.

{\it Case 2:} In this case we take  $u_1 \in   V_3^3(\ov{ M}_0)$ and assume that $ D_9 \cup D_8 \neq \emptyset$, where
$$
D_{9}:= \{ i\in I: i \in B_z  \mbox{ with } \# B_z =1\};  \qquad D_{8}:=  I \cap \mathcal{I}_{in}.
$$
Here  we define the following vector field:
  $$W_4^2:= W_4^0 + (1/M_2^2) \sum_{i\in D_8 \cup D_9} \l_i \frac{\partial \varphi_i}{\partial \l_i} .$$
We point out  that, this pseudogradient increases the $\mu_i$ for $i\in D_8 \cup D_9$, but does not increase the  $\mu_{\max}:= \mu_{p+q}$  since $p+q \notin I$.  Furthermore observe that
$$  \langle - \n J_K(u), \sum_{i\in D_8 \cup D_9} \l_i \frac{\partial \varphi_i}{\partial \l_i} \rangle \geq c \sum_{i\in D_8 \cup D_9} \Big( \frac{1}{\mu_i} + O\Big( \sum_{j =1}^{p+q} \frac{1}{\l_j^3} + \sum_{\ell \notin I} \e_{i\ell} \Big) \Big).$$
Hence the result follows as the first case.

Next we set
$$D_{10}:= \{i\in I:  \sum_{k\in I; k\neq i} \e_{ki} \leq m_1 q /\l_i\} \neq \emptyset, \quad  \mbox{ where  } m_1 \mbox{  is a small constant}.$$

{\it Case 3:} In this case we take  $u_1 \in   V_3^3(\ov{ M}_0)$ and  assume that  $ D_9 \cup D_8 = \emptyset$. That is we have that $I\subset \mathcal{I}_b$ and that  $\#B_z \neq 1$ for each $z$  critical point of $K_1$. Furthermore  we assume that $D_{10} \neq \emptyset$.

Next we recall that in  this case, for each $z$ such that  $\#B_z \geq 2$, $z$ has to be a local maximum point with $\partial K/\partial \nu > 0$ (which implies that $\mu_i$ and $\l_i$ are of the same order). Hence one can use  the same pseudogradient defined in  Case 2 (by  replacing   $D_{8} \cup D_{9}$ by  $D_{10}$). Hence
for $i \in D_{10}$, using Proposition \ref{p:33}, we derive that
$$  \langle - \n J_K(u),  \l_i \frac{\partial \d_i}{\partial \l_i} \rangle \geq \frac{c}{\l_i} + O( \sum_{j\neq i} \e_{ij}) \geq \frac{c}{\l_i} + c \sum_{j\neq i; j\in I} \e_{ij} + O( \sum_{j\notin I} \e_{ij})$$
and the proof follows as the previous cases.

{\it Case 4:}  $u_1 \in   V_3^3(\ov{ M}_0)$ and  $I\subset I_b$, $\#B_z \neq 1$ for each $z$ and $D_{10} = \emptyset$.\\
In this case, for each $z$ such that  $\#B_z \geq 2$, $z$ has to be a local maximum point with $\partial K/\partial \nu > 0$ (which implies that $\mu_i$ and $\l_i$ are of the same order). Let $z_1,\cdots,z_\ell$ be such that $\# B_{z_j} \geq 2$. Thus, the function $u$ can be written as
$$ u := \sum_{j=1}^\ell u_j + u_{\ell +1} \quad \mbox{ where } \quad u_j := \sum_{i\in B_{z_j}} \a_i \varphi_i  \mbox{ for } j\leq \ell \quad \mbox{ and } \quad u_{\ell +1} := \sum_{i\notin I} \a_i \varphi_i .$$
Notice that, for $j \leq \ell$, it follows that $u_j \in V(z_j, \#B_{z_j}, \eta, \e, \ov{M}_0)$ and in Lemma \ref{simple}, we have constructed a pseudogradient $W(z_j, \#B_{z_j})$ in this region. Now, we define
\be\label{chv44}W_4^4:= W_4^0 + \frac{1}{M_2^2} \sum_{j=1}^\ell W(z_j,\#B_{z_j}) (u_j).\ee
Observe that, by Lemma \ref{simple}, we have
$$  \langle - \n J_K(u),   W(z_j,\#B_{z_j}) (u_j)  \rangle \geq c  \sum_{k\in B_{z_j}} \Big( \sum_{r\neq k; r\in B_{z_j}} \e_{kr} ^\frac{n-1}{n-2} + O\Big( \sum_{r\notin B_{z_j}, r\in \mathcal{I}_{in}} \e_{kr} + \sum_{ r\notin B_{z_j}, r\in \mathcal{I}_{b}} \frac{1}{\l_k} | \frac{\partial \e_{kr}}{\partial a_k} | \Big)\Big) .$$
Furthermore we notice  that, for $r\notin B_{z_j}$ and  $r\in \mathcal{I}_{in}$, $(i)$ either $r \geq i_0$ and therefore the $\e_{kr}$ exists already  in $\ov{\G}_4$ or  $(ii)$ $r < i_0$ and, using Lemma \ref{pointint}, it follows that $a_r$ is close to a critical point  $y$ of $K$ in $\mathbb{S}^n_+$ which implies that $\e_{kr} \leq c(1/\l_k^{n-2} + 1/\l_r^{n-2})$.
Next  for $r\notin B_{z_j}$ and $r\in \mathcal{I}_{b}$, three situations  may occur
\begin{itemize}
  \item[(i)]
  $r \geq j_0$ and therefore the $\e_{kr}$ exists already  in $\ov{\G}_4$.
  \item[(ii)]
  $r < j_0$ and $r \notin I$. In this case it follows that $ \e_{kr }  \leq M_2/\l_r$ and thus (since $\l_r \geq  M_0 ^{1/q+p - 1 } \l_k$ for each $k \in I$) we have that
\begin{align}
\frac{1}{\l _k} | \frac{\partial \e_{kr}}{\partial a_k} | & \leq {c}{ \l_r d(a_r,a_k) } \e_{kr } ^{\frac{n}{n-2}} \leq {c}\sqrt{\frac{ \l_r }{\l_k}} \e_{kr } ^{\frac{n-1}{n-2}}
\leq \frac{c \, M_2^{\frac{n-1}{n-2}}}{ \l_k ^{1/2}\l_r ^{1/2 + 1/(n-2)}} \nonumber\\
& \leq c \frac{M_2 ^{\frac{n-1}{n-2}} } { M_0 ^{(^{1/2 + 1/(n-2)})/q+p - 1 }}\frac{1}{\l_k ^{1 + 1/(n-2)}} = o \Big(\big( \frac{m_1}{\l_k}\big)^{(n-1)/(n-2)}\Big)\label{M0M2}
\end{align}
(by using \eqref{M0M2M4}).
  \item[(iii)]
$r < j_0$ and $r \in I$. In this case, it follows that $a_r \in B_{z_\ell}$ with $\ell \neq k$ and therefore we deduce that $ | a_k - a_r | \geq c > 0$. Hence we get  $$ \frac{1}{\l _k} | \frac{\partial \e_{kr}}{\partial a_k} |  = O\Big( \frac{1}{\l_k ^{n-1} }+ \frac{1}{\l_r ^{n-1} }\Big).$$
\end{itemize}
Using \eqref{Gamma4},\eqref{chv44}, the previous estimates and the fact that $D_{10} = \emptyset$, the lemma follows in this case.
\end{pf}

\begin{pfn}{\bf of Proposition \ref{pp:champ2} }
The required pseudogradient will be  a convex combination of the ones defined in the previous lemmas. Each one is bounded and satisfies Claim $(i)$. Furthermore, the only case where $\mu_{\max}$ increases is the region $\mathcal{W}$. Finally, Claim $(ii)$ follows from the first one and the estimate of $\| \ov{v} \| ^2$ which is small with respect to the lower bound of Claim $(i)$. Concerning the last claim, it follows easily from the definition of the pseudogradient. This achieves the proof of Proposition \ref{pp:champ2}.
\end{pfn}

\subsubsection{ Critical points at Infinity and their topological contribution }

For $\e_0$ a small number, we define the following neighborhood of the  cone  of positive solutions of the sphere in $H^1(\mathbb{S}^n_+)$:
$$
 V_{\e_0}(\Sigma^+):= \{  u \in \Sigma; \, J_K(u)^{(2n-2)/(n-2)} e^{2J(u)}|u^-|_{L^{{2n}/{(n-2)}}}^{{4}/{(n-2)}} \, < \, \e_0 \}, \quad \mbox{ where } u^-:= \max(0, -u).
 $$
This set is for $\e_0$ small enough invariant under the gradient flow lines of the Euler Lagrange functional $J_K$. Namely we prove that
\begin{lem}\label{l:51}
For $\e_0 > 0$ small enough, the set $V_{\e_0} (\Sig ^+)$ is invariant under the flow generated by $-\n J_K$.
\end{lem}
\begin{pf} We will write $J$ instead of $J_K$. For $w\in L^{2n/(n+2)} (\mathbb{S}^n_+)$, we denote by $\mathcal{L}^{-1} (w)$ the solution of the following PDE:
$$ \begin{cases} \mathcal{L}u := -\D u + \frac{n(n-2)}{4} u = w \quad \mbox{in } \mathbb{S}^n_+, \\
 \partial u / \partial \nu = 0 \quad  \mbox{on } \partial \mathbb{S}^n_+ . \end{cases}$$
 Furthermore, it holds

 \begin{align} &  | u | _{L^{2n/(n-2)} } \leq c \| u \| _{H^1} \leq c | w |_{L^{2n/(n+2)} } \nonumber \\
&  | \mathcal{L} ^{-1}(K|u |^{4/(n-2)}u) | _{L^{2n/(n-2)} } \leq c | u |_{L^{2n/(n-2)} } ^{(n+2)/(n-2)}  .\label{99} \end{align}
Suppose $u_0 \in V_{\e_0} (\Sig ^+)$ and consider
$$
\begin{cases}
\frac{du(s)}{ds}= -\n J(u(s)) = -2J(u)\biggr(u-J(u)^{n/(n-2)} \mathcal{L} ^{-1}(K|u|^{4/(n-2)}u)\biggr)\\
u(0)=u_0.
\end{cases}
$$
Then
$$ e^{2\int_0^s J(u(t)) dt}u(s)= u_0 +2\int_0^s e^{2\int_0^tJ(u(y)) dy } J(u(t))^{\frac{2n-2}{n-2}}\mathcal{L} ^{-1}(K|u(t)|^{4/(n-2)}u(t))dt,
$$
$$ u^-(s) \leq e^{-2\int_0^s J(u(t)) dt}\biggr( u_0^- +2\int_0^se^{2\int_0^tJ(u(y)) dy} J(u(t))^{\frac{2n-2}{n-2}}\mathcal{L} ^{-1}(K(u^-(t))^{\frac{n+2}{n-2}})dt\biggr):= e^{-2\int_0^s J(u)}f(s).$$
 Setting
$$
F(s)=e^{-\frac{4n}{n-2}\int_0^s J(u(t)) dt } | f (s) |_{L^{{2n}/{(n-2)}}}^{{2n}/{(n-2)}}
\quad \mbox{ which implies that }\quad  |u^-(s)|_{L^{{2n}/{(n-2)}}}^{{2n}/{(n-2)}} \leq F(s).$$
Recall that, if $u_0^{-} =0$ then $u(s)$ is positive for all $s$. Hence, we can assume that $u_0^-\ne 0$ and we
want to prove that $F$ is a decreasing function. Observe that
\begin{align*}
F'(s)&= -\frac{4n}{n-2} J(u(s))e^{-\frac{4n}{n-2}\int_0^s J(u)}
|f(s) |_{L^{{2n}/{(n-2)}}}^{{2n}/{(n-2)}}+ e^{-\frac{4n}{n-2}
\int_0^s J(u)}\frac{2n}{n-2} \int_{\mathbb{S}^n_+} f' (s) f(s) ^{\frac{n+2}{n-2} } dx \\
&\leq \frac{2n}{n-2} e^{-\frac{4n}{n-2}\int_0^s J(u)}
\left[-2J(u(s))|u_0^-|_{L^{{2n}/{(n-2)}}}^{{2n}/{(n-2)} } +
\int_{\mathbb{S}^n_+}f'(s) f(s) ^{\frac{n+2}{n-2} } dx\right] \qquad (\mbox{using } f(s) \geq u_0^{-}).
\end{align*}
Notice that $f'(0) = u_0^{-}$ and therefore
\begin{eqnarray*}
\big|\int_{\mathbb{S}^n_+}f'(s)f(s)^{\frac{n+2}{n-2}}\big| dx \leq c\int_{\mathbb{S}^n_+} | f'(s) |
|u_0^-|^{\frac{n+2}{n-2}}+c\int_{\mathbb{S}^n_+} | f'(s) | \Big(\int_0^s |f'(t)|dt \Big)^{\frac{n+2}{n-2}} dx.
\end{eqnarray*}
But, we have (using \eqref{99})
\begin{align*}
\int_{\mathbb{S}^n_+}(u_0^-)^{\frac{n+2}{n-2}}  | f'(s) | dx & = \int_{\mathbb{S}^n_+}(u_0^-)^{\frac{n+2}{n-2}}\Big( 2 e^{2\int_0^sJ(u)}
J(u(s))^{\frac{2n-2}{n-2}}\mathcal{L} ^{-1}(K(u^-(s))^{\frac{n+2}{n-2}})\Big) dx\\
& \leq CJ(u(s))^{\frac{2n-2}{n-2}} e^{2\int_0^sJ(u)} |u_0^-|_{L^{{{2n}/{(n-2)}}}}^{{(n+2)}/{(n-2)}} |u^-(s)|_{L^{{{2n}/{(n-2)}}}}^{{(n+2)}/{(n-2)}},
\end{align*}
and we also have (using the fact that $J(u(s))$ is a decreasing function)
\begin{align*}
\int_{ \mathbb{S}^n_+ }|f'(s)| & \big( \int_0^s|f'(t))|dt \big)^{\frac{n+2}{n-2} }dx  \leq cs^{\frac{4}{n-2}}
\int_{\mathbb{S}^n_+}|f'(s)|\int_0^s|f'(t)|^{\frac{n+2}{n-2}}dt dx \\
 & \leq c s^{\frac{4}{n-2}} e^{\frac{4n}{n-2}sJ(u_0)} J(u_0)^{\frac{2n-2}{n-2}\frac{2n}{n-2}} |u^-(s)|_{L^{{{2n}/{(n-2)}}}}^{{(n+2)}/{(n-2)}}  \int_0^s  |u^-(t)|_{L^{{{2n}/{(n-2)}}}}^{(n+2)^2/(n-2)^2} dt.
\end{align*}
Hence, if $|u^-(s)|_{L^{{2n}/{(n-2)}}} \leq 5|u_0^-|_{L^{{2n}/{(n-2)} }}$,
for $0\leq s\leq 1$, we derive that
\begin{align*}
 F'(s) \leq \frac{4 n}{n-2} e^{-\frac{4n}{n-2}\int_0^s J(u)}
|u_0^-|_{L^{{2n}/{(n-2)}}}^{{2n}/{(n-2)} } \Big( - J(u(s)) + c \, J(u_0)^{\frac{2n-2}{n-2}} e^{2 J(u_0) } | u^- _0 | _{L^{{2n}/{(n-2)}}}^{{4}/{(n-2)}}  \\
+ c\, \Big( J(u_0)^{\frac{2n-2}{n-2 } } e^{2 J(u_0)} | u^- _0 | _{L^{{2n}/{(n-2)}}}^{4/{(n-2)}} \Big)^{2n/(n-2)}\Big)
\end{align*}
Finally, since $\inf J > c > 0$, using the fact that $u_0\in V_{\e_0}(\Sig^+)$, that is, $J(u_0)^{\frac{2n-2}{n-2}} e^{2J(u_0)}|u_o^-|_{L^{{2n}/{(n-2)}}}^{{4}/{(n-2)}} <\e_0 $,
and $\eta$ is small enough, then $F'(s) \leq 0$, for $0\leq s\leq 1$. Therefore
$J(u(s))^{\frac{2n-2}{n-2}} e^{2J(u(s))}|u(s)^-|_{L^{{2n}/{(n-2)}}}^{{4}/{(n-2)}} <\e_0 $,
and our result follows.
\end{pf}



Next using a partition of the unity, one can define the vector field $W$ of Proposition \ref{pp:champ2} globally by gluing it to the negative gradient $- \n J$ outside the $V(q,p,m,\e)$'s. Let us denote the resulting global vector field by $Y$ and define a new vector field by setting:
$$
X(u):= \, Y(u) \, -  \, <Y(u),u> u \qquad  \mbox{ for } u \in V_{\e_0}(\Sigma^+).
$$
We then have
\begin{cor}\label{c:pseudo}
Assume that  $K$ satisfies $(H1)$, $(H2)$ and $(H3)$. Then using  Propositions \ref{pp:champ2m1}, \ref{pp:champ2} and arguing as in the above Lemma, one  proves that   for $\e_0$ small enough, $X$ is a pseudogradient of $J$ which preserves $ V_{\e_0}(\Sigma^+)$. Moreover  the \emph{critical points at infinity} of $X$ lie in subsets $\mathcal{W}$ (see the formula \eqref{setW} for a definition)
\end{cor}
Next we perform a Morse type reduction in the subsets $\mathcal{W}$. Namely we prove

\begin{lem}\label{l:morse}
For $u = \sum_{i=1}^{q} \a_i \d_{a_i,\l_i} + \sum_{q+1}^{p+q} \a_i \varphi_{a_i,\l_i}\in \mathcal{W}$, we define
$$
D_4:= \{ i \leq q : a_i  \mbox{ is close to } z \mbox{ with } \frac{\partial K}{\partial \nu}(z) = 0\} \, \, \& \, \, D_5 := \{ i \leq q : a_i  \mbox{ is close to } z \mbox{ with } \frac{\partial K}{\partial \nu}(z) > 0\}.
$$
Then the functional $J_K$ expands as follows

\begin{align}
J_K(u) & = \frac{(\sum_{i\leq q }\a _i^2+2 \sum_{i > q}\a _i^2)S_n^{2/n}}{(\sum_{i\leq q }\a _i^{\frac{2n}{n-2}}K(a_i)+2\sum_{i > q } \a _i^{\frac{2n}{n-2}}K(a_i))^{\frac{n-2}{n}}}\Big( 1- c \sum_{i > q } \frac{\D K(y_i)}{\l_i^2}  + c  \sum_{i \in D_5} \frac{1}{\l _i}\frac{\partial K}{\partial  \nu}(z_i) \nonumber\\
 & + c  \sum_{i \in D_4} \Big( \frac{c_7}{\l _i}\frac{\partial K}{\partial  \nu}(a_i) - c_6 \frac{\D K(a_i)}{\l_i^2}\Big)  + o \Big( \sum_{i \in D_5} \frac{1}{\l _i}+ \sum_{i \in D4}  \frac{1}{\l_i^2}+ \sum_{i > q } \frac{1}{\l_i^2} \Big) \Big)\nonumber \\
 & = S_n^{2/n} \Big( \sum_{i \leq q } \frac{1}{K(z_i)^{\frac{n-2}{2}}} + 2 \sum_{i > q } \frac{1}{K(y_i)^{\frac{n-2}{2}}} \Big)^\frac{2}{n}
 \Big( 1 - \| \a\|^2 + \sum _{i=1}^{p+q} \big( | A_i^-|^2 - | A_i^+| \big) - c \sum_{i > q } \frac{\D K(y_i)}{\l_i^2}  \label{expJ3}\\
 &  + c    \sum_{i \in D_5} \frac{1}{\l _i}\frac{\partial K}{\partial  \nu}(z_i) + c  \sum_{i \in D_4} \Big( \frac{c_7}{\l _i}\frac{\partial K}{\partial  \nu}(a_i) - c_6 \frac{\D K(z_i)}{\l_i^2}\Big) + o \Big(  \sum_{i \in D_5} \frac{1}{\l _i}+ \sum_{[i > q]\cup [i \in D_4] } \frac{1}{\l_i^2} \Big) \Big),  \nonumber
\end{align}
where $S_n$ is defined in Proposition \ref{p:311} (it represents the level of one boundary bubble),  $\a \in \R^{q+p-1} $, $(A_i^+,A_i^-)$ are the local coordinates of the parameters $(\a_1,\cdots,\a_{p+q})$ and $a_i$. This expansion will be called the {\it Morse Lemma at Infinity} of $J_K$ near its critical point  at infinity. Note that we loose an index for the parameter $\a$ since the functional $J_K$ is homogenous with respect to this parameter.
\end{lem}


From Propositions \ref{pp:champ2m1}, \ref{pp:champ2} and Lemma \ref{l:morse}, we derive the characterization of \emph{critical points at infinity} and identify their level sets. Namely we have:
\begin{cor}\label{caract2}
Assume that  $K$ satisfies $(H1)$, $(H2)$ and $(H3)$. Then, in $V(m,q,p,\e)$, the critical points  at infinity of $J_K$ are in one to one correspondence with the collections of $q$ critical points  $z_\ell$'s of $K_1$ satisfying: either $z_\ell$ is a local maximum point with $\partial K/ \partial \nu (z_\ell ) >0$ or $\partial K/ \partial \nu (z_\ell )=0$ and $\D K(z_\ell ) < 0$ and $p$ critical points  $y_r$'s of $K$ in $\mathbb{S}^n_+$ satisfying $\D K(y_r) < 0$. We will denote such a critical point  at infinity  by $(z_1, \cdots,z_q, y_{q+1}, \cdots,y_{q+p})_\infty$. Such a critical point  at infinity is at the level (see \eqref{expJ3})
$$C_\infty (z_1, \cdots,z_q, y_{q+1}, \cdots,y_{q+p}):= S_n^{2/n} \Big( \sum_{i=1}^q \frac{1}{K(z_i)^{(n-2)/2}} + \sum_{i=q+1}^{q+p} \frac{2}{K(y_i)^{(n-2)/2}}\Big)^{2/n}. $$
In particular, it holds that
$$C_{\min}^{(2p + q),\infty} :=\frac{\big((2p + q ) S_n\big)^{2/n}}{K_{\max}^{(n-2)/n}} \leq C_\infty (z_1, \cdots,z_q, y_{q+1}, \cdots,y_{q+p}) \leq  \frac{\big( (2p + q ) S_n\big)^{2/n} }{K_{\min}^{(n-2)/n}}:= C_{\max}^{(2p + q),\infty}  $$
Furthermore, for such a critical point  at infinity, we associate an index (which corresponds to the number of the decreasing directions for $J_K$ by using the Morse Lemma at infinity, see \eqref{expJ3}) $$\ i_\infty (z_1, \cdots,z_q, y_{q+1}, \cdots,y_{q+p}) := q+p-1 + \sum_{i=1}^q (n-1 - morse(K_1, z_i)) + \sum_{i= q+1} ^{q+p} (n- morse (K,y_i)). $$
Such an index will be called the {$i_{\infty}$-index of such a critical point  at infinity}.
\end{cor}

Next as consequence of the above corollary and the Morse reduction in Lemma \ref{l:morse} we compute the topological contribution of the \emph{critical points at infinity} to the difference of topology between the level sets of the functional $J_K$. Namely we have

\begin{lem}\label{l:difftop}
Let $\tau_{\infty}$ be a critical point at infinity at the level  $C_{\infty}(\tau_{\infty})$ with index  $i_{\infty}(\tau_{\infty})$. Then for $\theta$ a small positive number and  a field $\mathbb{F}$, we have that

$$
 H_l(J_K^{C_{\infty}(\tau_{\infty}) + \theta},   J_K^{C_{\infty}(\tau_{\infty}) - \theta }; \mathbb{F} ) =
\begin{cases}
   \mathbb{F} & \mbox{ if  } \quad  l =  i_{\infty}(\tau_{\infty}), \\
  0, & \mbox{otherwise}.
\end{cases}
$$
where $H_l$ denotes the $l-$dimensional homology group with coefficient in the field $\mathbb{F}$.

\end{lem}

\section{Proof of the main results}

This section is devoted to the proof of Theorems \ref{t:pinching1}, \ref{t:pinching2} and \ref{t:pinching3}. The proof of these theorems is based  on the characterization of the critical points at infinity in   Corollary \ref{caract2} and the computation of their contribution to the difference of topology in Lemma \ref{l:difftop}. It also uses two  deformation lemmas. The first one is an abstract lemma, which is inspired by Proposition 3.1 in \cite{MM19}. It reads as follows:
\begin{lem}\label{deform}
Let $\underline{A}$ and $\ov{A}:=  ({ K_{\max}}/{ K_{\min}})^{(n-2)/n} \, \underline{A}$. Assume that $J_K$ does not have any critical point  nor critical point  at infinity in the set $J_K^{\ov{A}} \setminus J_K^{\underline{A}}$ where $J_K^A := \{ u: J_K(u) < A\}$. Then for each $c\in [ \underline{A},\ov{A}]$, the level set  $J_K^c$ is contractible.
\end{lem}
\begin{pf}
First, since we assumed that $J_K$ does not have any critical point  nor critical point  at infinity in $\Sigma^+$ between the levels $\ov{A}$ and $\underline{A}$, we have  that $J_K^{\ov{A}}$ retracts by deformation onto $J_K^{\underline{A}}$. Indeed such a retraction can be realized by following the flow lines of a decreasing pseudogradient  $Z_K$ for $J_K$. Let $\phi_K$ denote  the one parameter group corresponding to this pseudogradient.  For each $u\in \Sigma^+$, we denote by $s_K(u)$ the first time such that $\phi_K(s_K(u), u ) \in J_K^{\underline{A}}$.\\
Secondly we recall that, for $K\equiv 1$, the only critical points  of $J_1$ are minima and lie  in the bottom level $S_n$. Furthermore, for each $A > S_n$, the set $J_1^A$ is a contractible one. Indeed by following the flow lines of  a decreasing pseudogradient $Z_1$ of the Yamabe functional $J_1$, each flow line, starting from $u\in\Sigma^+$, will reach the bottom level $S_n$.
Let us denote by $\phi_1$ the one parameter group corresponding to $Z_1$.\\
Next we notice  that, we have
$$ (1/ K_{\max}^{(n-2)/n}) J_1 (u) \leq J_K (u) \leq   (1/ K_{\min}^{(n-2)/n}) J_1 (u)  \quad \mbox{ for each } u\in \Sigma, $$
which implies that
$$ J_K ^{\underline{A}} \subset  J_1 ^{A'} \subset  J_K ^{\ov{A}} \quad \mbox{ where  } A':=  K_{\max}^{(n-2)/n} \underline{A}.$$
Furthermore we observe that  for each $u \in \Sigma^+$, there exists a unique $s_1(u)$ satisfying $\phi_1(s_1(u), u ) \in J_1^{A'}$.\\
Next we define the following map:
$$ F:= [ 0,1] \times J_1^{A'} \to  J_1^{A'} \, ; \quad F(t,u) := \phi_1(s_1(\phi_K(t\, s_K(u), u )), \phi_K(t\, s_K(u), u ) ).$$
We notice that $F$ is well defined and continuous and satisfies the following properties:
\begin{itemize}
\item For $t=0$, we have $\phi_K(0, u )=u$. Furthermore, for each $u\in  J_1^{A'} $, we have $s_1(u) = 0$. Therefore, for each $u\in  J_1^{A'} $, we get $F(0,u)= \phi_1(0,u) = u$.
\item For $t=1$, we have $\phi_K( s_K(u), u ) \in J_K^{\underline{A}} \subset J_1^{A'}$ (by the definition of $s_K$) which implies that $s_1(\phi_K( s_K(u), u ))= 0$ and therefore $F(1,u)= \phi_1(0,\phi_K( s_K(u), u )) = \phi_K( s_K(u), u ) \in J_K^{\underline{A}}$ for each $u \in J_1^{A'}$.
\item If $u\in J_K^{\underline{A}}$, then $s_K(u)=0$ which implies that $\phi_K(t\, s_K(u), u ) = \phi_K( 0, u ) =u$. Therefore $F(t,u) = \phi_1(s_1(u),u) = \phi_1(0,u)=u$ for each $u\in J_K^{\underline{A}}$ and each $t\in [0,1]$ (we used $s_1(u)=0$ since $u\in J_K^{\underline{A}}\subset J_1^{A'}$).
\end{itemize}
Thus  $ J_1^{A'}$ retracts by deformation onto $J_K^{\underline{A}}$, a fact which  provides the claim of the lemma  since $J_1^{A'}$  itself is a contractible set.
\end{pf}

The second deformation lemma is a consequence of the previous one, the assumptions  $(H1)$, $(H2)$, $(H3)$ of this paper and an appropriate pinching condition for the function $K$. To state it we set the following notation:
$$
\mbox{ for } \ell \in \N, \quad C^{\ell, \infty }_{\max} := (\ell S_n)^{2/n} / K_{\min} ^{(n-2)/n} \quad \&  \quad C^{\ell, \infty }_{\min} := (\ell S_n)^{2/n} / K_{\max} ^{(n-2)/n}.
$$
We recall that it follows from Corollary \ref{caract2} that  the level of critical points at Infinity corresponding to $q$ boundary points and $p$ interior points such that $ q + 2p= \ell$  lie between $ C^{\ell, \infty }_{\min}$ and $  C^{\ell, \infty }_{\max}$. \\
Our second deformation lemma reads as follows:

\begin{pro}\label{deform2}

For $k\in \N$ a fixed integer, let $ 0 < K \in C^3( \ov{\mathbb{S}^n_+})$ satisfying the conditions  $(H1)$, $(H2)$, $(H3)$ and the pinching condition $K_{\max} / K_{\min} < ((k+1)/k)^{1/(n-2)}$.\\
Assume that $J_K$ does not have any critical point  under the level $C^{k+1,\infty} _{\min}$. Then, for every  $1\leq \ell \leq k$ and every  $c\in (C^{\ell,\infty}_{\max}, C^{\ell +1,\infty}_{\min} )$, the sublevel $J_K^c$ is a contractible set.
\end{pro}
\begin{pf}
Since we assumed that $K_{\max} / K_{\min} < ((k+1)/k)^{1/(n-2)}$, it follows that, for each $1\leq \ell \leq k$, we have $(k+1)/k \leq (\ell +1)/\ell$ and
$$C^{\ell , \infty }_{\max} <  C^{\ell , \infty }_{\max} (K_{\max} / K_{\min} )^{(n-2)/n} <  C^{\ell +1 , \infty }_{\min} .$$
The proof follows then from Lemma  \ref{deform} by taking $\underline{A}= C^{\ell, \infty }_{\max} + \g$ with a small $\g >0$ so that
$\ov{A} < C^{\ell +1, \infty }_{\min} $. Indeed  between the levels $\underline{A}$ and $ \ov{A}$ the functional $J_K$ does not have any critical point  nor critical point  at infinity.
\end{pf}

\bigskip
\noindent
Next we start the proof of our existence results by proving   Theorem \ref{t:pinching2}.

\begin{pfn}{\bf of Theorem \ref{t:pinching2}}
Arguing by contradiction we assume that the functional $J_K$ does not have any critical point under the level $C_{\min}^{2,\infty}$. Hence it follows  from   Proposition \ref{deform2} (with $k=1$) that  under the assumption of Theorem \ref{t:pinching2}, we have  that $J^{C_{\max}^{1,\infty} +\g}$ is a contractible set, for $\g$  a small constant. Moreover it  is a retract by deformation of   $C_{\min}^{2,\infty}$. Furthermore  follows from corollary \ref{caract2} that critical points at infinity under the level $C_{\min}^{2,\infty}$ are in one to one correspondence with critical points of $K_1$ in  $\mathcal{K}^{+}_b\cup \mathcal{K}^{0,-}_b$. Then it follows from Lemma \ref{l:difftop} and the Euler-Poincar\'e theorem that:
$$1 = \chi( J^{C_{\min}^{2,\infty} +\g}) = \sum_{ z\in \mathcal{K}^{+}_b\cup \mathcal{K}^{0,-}_b} (-1)^{n-1- morse (K_1,z)}$$ which contradicts  the  assumption $(b)$ of Theorem \ref{t:pinching2}. Hence the existence of at least one critical point  of $J_K$.
 \end{pfn}

 \begin{pfn}{\bf of Theorem \ref{t:pinching3}}
Assuming that $J_K$ does not have any critical point under the level $C_{\min}^{3,\infty}$,  we derive, using Proposition \ref{deform2} (with $k=2$), the level sets  $J_K^{C_{\max}^{1,\infty} +\g}$ and $J_K^{C_{\max}^{2,\infty} +\g}$ are contractible sets. Then it follows from the properties of the Euler-Characteristic, see Proposition 5.7, pp.105 in \cite{Dold},  that
$$
1 = \chi(J_K^{C^{2,\infty}_{\max}+\g}) \, = \, \chi(J_K^{C^{2,\infty}_{\max}+\g}, J_K^{C^{1,\infty}_{\max}+\g}) \, + \, \chi(J_K^{C^{1,\infty}_{\max}+\g}).
$$
That is   $  \chi(J_K^{C^{2,\infty}_{\max}+\g}, J_K^{C^{1,\infty}_{\max}+\g}) \, = \, 0.$
Moreover it follows from Corollary  \ref{caract2} that the critical points at infinity between these two levels   are $(z_i,z_j)_\infty$ with $z_i \neq z_j \in \mathcal{K}^{+}_b\cup \mathcal{K}^{0,-}_b$ and $y_\infty$ with $y \in \mathcal{K}^{-}_{in}$. Thus, it follows from Lemma \ref{l:difftop} and the Euler-Poincar\'e theorem that
$$\sum_{z_i \neq z_j \in \mathcal{K}^{+}_b\cup \mathcal{K}^{0,-}_b} (-1)^{1+ \i(z_i) + \i(z_j)} + \sum_{y \in \mathcal{K}^{-}_{in}} (-1)^{\i(y)} = 0$$
where $\i(z_k):= n-1-\mbox{morse}(K_1,z_k)$ and $\i(y):= n-\mbox{morse}(K,y)$. \\
Observe that, the first term is exactly $-A_2$ defined in Lemma \ref{indicesz}. Hence, the previous equality contradicts the assumption $(ii)$ of the theorem. The proof is thereby completed.
\end{pfn}

\begin{pfn}{\bf of Theorem \ref{t:pinching1}}
We first observe that, under the assumption of the theorem,   if  $A_1 \neq  1$ or   respectively $A_1 =1$ and $ B_1 \neq -k$, where $\#( \mathcal{K}^{+}_b\cup \mathcal{K}^{0,-}_b)= 2k+1$,  the existence of  at least one solution to Problem  $(\mathcal{P})$ follows from Theorem \ref{t:pinching3},   respectively Theorem \ref{t:pinching2}. Hence we will assume that $A_1=1$ and $B_1=-k$ and notice that
$$
\# ( \mathcal{K}_{in}^{-}) \, = 2 r + k, \mbox{ where } r \in \mathbb{N}_0,
$$
and  there are $r$ even numbers $\i(y_j)$'s and $r+k$ odd numbers $\i(y_j)$'s. \\
Next arguing as in the proof of Theorem \ref{t:pinching3} using the assumption on $K_{\max}/K_{\min}$ and Proposition \ref{deform2}, we deduce that $J_K^{C_{\max}^{3,\infty} +\g}$ and $J_K^{C_{\max}^{4,\infty} +\g}$ are contractible sets.
 Using Corollary \ref{caract2}, we derive that the critical points at infinity  whose level are lying between these values   are :
\begin{itemize}
\item $(z_i,z_j,z_r,z_t)_\infty$ with different $z_i$'s which belong to $\mathcal{K}^{+}_b\cup \mathcal{K}^{0,-}_b$,
\item  $(z_i,z_j,y)_\infty$ with $y \in \mathcal{K}^{-}_{in}$ and $z_i \neq z_j \in \mathcal{K}^{+}_b\cup \mathcal{K}^{0,-}_b$,
\item $(y_i,y_j)_\infty$ with $y_i \neq y_j \in \mathcal{K}^{-}_{in}$.
\end{itemize}
Hence arguing as above  we derive that

\begin{align*}
 \sum_{z_i \neq z_j\neq z_r\neq z_t \in \mathcal{K}^{+}_b\cup \mathcal{K}^{0,-}_b} & (-1)^{3+ \i(z_i) + \i(z_j) + \i(z_r) + \i(z_t)} \\
&  + \sum_{y \in \mathcal{K}^{-}_{in};z_i \neq z_j \in \mathcal{K}^{+}_b\cup \mathcal{K}^{0,-}_b} (-1)^{2 + \i(z_i) + \i(z_j) +\i(y)} + \sum_{y_i \neq y_j \in \mathcal{K}^{-}_{in}} (-1)^{1+ \i(y_i)+\i(y_j)} = 0.\end{align*}
Observe that, the first term is exactly $-A_4$, the second one is $A_2\times B_1$ and the third one is $-B_2$ (defined in Lemmas \ref{indicesz} and \ref{indicesy}). Using the values of these terms (given in Lemmas \ref{indicesz} and \ref{indicesy}), we obtain that $$r+k=0$$ which implies that $r=k=0$. Now, from $r=k=0$, we get $\#(\mathcal{K}^{+}_b\cup \mathcal{K}^{0,-}_b) = 1$ and $\# \mathcal{K}^{-}_{in} = 0$. This leads to a  contradiction with the assumption that  $\#(\mathcal{K}^{+}_b\cup \mathcal{K}^{0,-}_b\cup \mathcal{K}^{-}_{in}) \geq 2$. Thereby  the proof of the theorem is completed.
\end{pfn}



\section{Appendix}

\subsection{Bubble estimates}
\begin{lem}\label{lem:1}
For $a\in \partial \mathbb{S}_+^n$, we have $\partial \delta_{a,\l} /\partial \nu = 0$ and therefore $\varphi_{a,\l} = \d_{a,\l}$. For $a\notin \partial \mathbb{S}_+^n$, we have
$$(i) \quad \d_{a,\l} \leq  \varphi_{a,\l} \leq 2 \d_{a,\l}\, \, ; \quad | \l \partial \varphi_{a,\l} /\partial \l | \leq c \d_{a,\l} \, \, ; \quad  | (1/\l) \partial \varphi_{a,\l} /\partial a^k | \leq c \d_{a,\l},$$
where $a^k$ denotes the $k$-th component of $a$.
$$ (ii) \qquad  \varphi_{a,\l}  = \d_{a,\l} + c_0 \frac{H(a,.)}{\l^{(n-2)/2}} + f_{a,\l}\qquad \mbox{ where }$$
$$ | f_{a,\l} | _{\infty} \leq \frac{c}{\l^{(n+2)/2} d_a^{n}} \, \, ; \quad | \l \frac{\partial f_{a,\l} }{\partial \l}  | _{\infty} \leq \frac{c}{\l^{(n+2)/2} d_a^{n}} \quad\mbox{ and } \quad | \frac{1}{\l} \frac{\partial f_{a,\l} }{\partial a^k}  | _{\infty} \leq \frac{c}{\l^{(n+4)/2} d_a^{n+1}},$$ where $d_a:= d(a,\partial \mathbb{S}^n_+$).
\end{lem}
\begin{pf}
Using a stereographic projection, we are led to prove the corresponding estimates on $\R^n_+$. We still denote by $G$ and $H$ the Green’s function and its regular part of Laplacian on $\R^n_+$ under Neumann boundary conditions. In this case, we have
$$\d_{a,\l}(x):= c_0 \frac{\l^{(n-2)/2}}{ (1+\l^2| x-a|^2)^{(n-2)/2}} \qquad \mbox{ and } \qquad H(a,x) := \frac{1}{ | x - \ov{a} | ^{n-2}}, $$ where $\ov{a}$ denotes the symmetric point  of $a$ with respect to $\partial \R^n_+$. Let $\psi:= \d_{a,\l} + \d_{\ov{a},\l}$. Easy computation implies that $\partial \psi /\partial \nu = 0$.\\
To prove the first inequality, let us consider $h:=  \varphi_{a,\l}  - \d_{a,\l} $. Hence we get $\D h = 0$ and $\partial h /\partial \nu = -\partial \d_{a,\l} /\partial \nu > 0$. Hence, using the Green's representation, we derive that $h > 0$ in $\R^n_+$.\\
For the second inequality, let us consider $h := \psi - \varphi_{a,\l}$. Easy computations imply that $\partial h /\partial \nu = 0$ and $-\D h = - \D \d_{\ov{a},\l} > 0$. Hence, $h >0$ in $\R^n_+$. The inequality follows from the fact that $\d_{\ov{a},\l} \leq \d_{{a},\l}$ in $\R^n_+$.\\
For the third one, let $g:= \l \partial \varphi_{a,\l} /\partial \l$, observe that $\partial g /\partial \nu = 0$ and $| \D g | \leq  ((n+2)/2)  \d_{a,\l}^{(n+2)/(n-2)}$. Now let us consider $h:=  ((n+2)/2) \psi \pm g$. It follows that $-\D h > 0$ and $\partial h /\partial \nu = 0$. Hence $h >0$ in $\R^n_+$ which gives the proof of the third inequality. The fourth one follows by the same way.\\
Concerning the second claim, it is easy to see that $\D  f_{a,\l} = 0$ and
\begin{align*}
 \frac{\partial  f_{a,\l}}{\partial \nu } & = - \frac{\partial  \d_{a,\l}}{\partial \nu } -  \frac{c_0}{\l^{(n-2)/2}} \frac{\partial  H(a,.)}{\partial \nu }= c_0 (n-2) \frac{\l^{(n+2)/2} d_a}{(1+ \l^2 | x-a|^2)^{n/2}} -  \frac{c_0}{\l^{(n-2)/2}} \frac{(n-2) d_a}{ | x-\ov{a} |^n}\\
 & = O \Big( \frac{d_a}{ \l^{(n+2)/2}  | x-a|^{n+2} } \Big).
 \end{align*}
Now, using the Green's representation, we get
\begin{align*}
| f_{a,\l} (x)| & \leq c \int_{\partial \R^n_+} G(x,y) |  \frac{\partial  f_{a,\l}}{\partial \nu } (y)| dy \leq \frac{c\, d_a}{ \l^{(n+2)/2}} \int_{\partial \R^n_+} G(x,y) \frac{1}{  | y-\ov{a} |^{n+2} }dy \\
& \leq \frac{c}{ \l^{(n+2)/2}d_a } \int_{\partial \R^n_+} G(x,y) \frac{1}{  | y-\ov{a} |^{n} }dy \leq \frac{c}{ \l^{(n+2)/2}d_a } \frac{H(a,x)}{d_a} \leq \frac{c}{ \l^{(n+2)/2}d_a^{n} }.
\end{align*}
This gives the first claim in $(ii)$. The other ones can be done by the same way.
\end{pf}


\begin{lem}\label{epsij1} 1) For each $i\neq j$, we have $$ - \l_i \frac{\partial \e_{ij}}{\partial \l_i} -  \l_j \frac{\partial \e_{ij}}{\partial \l_j } \geq 0\, \, \mbox{ and } - \l_i \frac{\partial \e_{ij}}{\partial \l_i}\geq c \e_{ij} \mbox{ if } \l_i \geq c \l_j \mbox{ or } \l_i d(a_i, a_j ) \geq 2. $$
2) Let $i,j \in I_b:=\{k:a_k\in \partial \mathbb{S}^n_+\}$ and let $\mu_i$ and $\mu_j$ be defined by \eqref{mui}.  Assume that $\mu_j \leq c' \mu_i$ for some constant $c'$, then: $(i)$ either there exists a constant $c''$ such that $\l_j \leq c'' \l_i$, $(ii)$ or $\l_i d(a_i,a_j ) \geq 2$.
\end{lem}

\begin{pf} The proof of the first assertion follows immediately from the definition of $\e_{ij}$. Concerning the second one, observe that, if $| \n K(a_i) | \geq c $ and $| \n K(a_j) | \geq c $, then it follows that $\mu_k$ and $\l_k$ are of the same order (that is: the ratio is bounded from above and below) for $k=i,j$. Hence the result follows in this case. In the other case, there exists $k\in \{i,j\}$ such that $a_k$ is close to a critical point $z$ of $K$ in $\partial \mathbb{S}^n_+$ (i.e. $\partial K/\partial \nu (z) =0$). Arguing by contradiction, assume that $\l_i d(a_i,a_j ) \leq 2$ and $\l_j / \l_i $ is very large. It follows that $a_i$ and $a_j$ are close to the same critical point $z$. Now we claim that:\\
{\bf Claim 1}: $\l_j | \n K(a_j) | $ is very large.\\
In fact, if it is not, we derive that $| \n K(a_j) | /\l_j \leq c/\l_j^2$ which implies that $1/\mu_j \leq c/\l_j^2$ and therefore $1/\mu_j$ is very small with respect to $1/\l_i^2 \leq 1/\mu_i$. This gives a contradiction and therefore our claim follows.\\
Since $z$ is a non degenerate critical point  of $K_1$, it follows that $\l_j d( a_j, z )$ is very large. Moreover, Claim 1 implies that $| \n K(a_j) | /\l_j \leq 1/\mu_j \leq c | \n K(a_j) | /\l_j $. Now we claim that:\\
{\bf Claim 2}: $\l_i | \n K(a_i) | \geq 1$ cannot occur.\\
To prove this claim, we assume that the inequality is true. Then we derive that $| \n K(a_i) | /\l_i \leq 1/\mu_i \leq 2 | \n K(a_i) | /\l_i$. Since $\mu_j \leq c' \mu_i$, we derive that $| \n K(a_i) | /\l_i  \leq c | \n K(a_j) | /\l_j$ and therefore $\l_j d( a_i, z ) \leq c  \l_i d( a_j, z )$ which implies that $d( a_i,z )  $ is very small with respect to $d(a_i,z) $ and therefore $d( a_i, z )$ is very small with respect to $ d( a_i, a_j ) $. Now observe that, since we assumed that  $\l_i | \n K(a_i) | \geq 1$,  it follows that $\l_i d( a_i,z ) \geq c$ and therefore $\l_i  d( a_i, a_j ) $ becomes very large
 which gives a contradiction. Hence Claim 2 follows.\\
Finally, we claim that \\
{\bf Claim 3}: $\l_i | \n K(a_i) | \leq 1$ cannot occur.\\
Arguing by contradiction we assume that  $\l_i d( a_i, z ) \leq c$. From $\mu_j \leq c' \mu_i$, we derive that $1/\l_i^2 \leq c| \n K(a_j) | /\l_j \leq c
d(a_j,z ) /\l_j$ and therefore $\l_j/\l_i \leq c \l_i  d(a_j, z )$, that is $\l_i  d( a_j, z )$ is very large. But we have $\l_i  d( a_j , a_i ) \leq 2$ and $\l_i  d( a_i, z ) \leq c$ which imply that $\l_i  d( a_j , z )$ is bounded. Hence we get a contradiction which completes the proof of Claim 3. \\
Hence the lemma is fully proven.
\end{pf}




\subsection{ Asymptotic expansion of the functional and its  gradient }

\begin{pro}\label{p:311}
Let $n\geq 5$ and $ u=\sum_{i\leq q} \a_i \d_i + \sum_{i > q} \a _i \varphi _i \in V(m,q,p,\e )$ be such that: $(i)$ $d( a_i, a_j ) \geq c > 0$ for every  $i\neq j$, $(ii)$ for $i > q$, $a_i$ is close to a critical point  $y_{j_i}$ of $K$ in $\mathbb{S}^n_+$ and $(iii)$ for $i \leq q$, $a_i$ is close to a critical point  $z_{j_i}$ of $K_1$ in $\partial \mathbb{S}^n_+$. Then the following expansion holds
\begin{align*}
J_K(u) & = \frac{(\sum_{i\leq q }\a _i^2+2 \sum_{i > q}\a _i^2)S_n^{2/n}}{(\sum_{i\leq q }\a _i^{{2n}/{n-2}}K(a_i)+2\sum_{i > q } \a _i^{{2n}/{n-2}}K(a_i))^{{n-2}/{n}}}\Big( 1 - 2 c_6 J_K(u)^{\frac{n}{n-2}} \sum_{i > q } \a_i ^\frac{2n}{n-2} \frac{\D K(a_i)}{\l_i^2} \\
& + J_K(u)^{\frac{n}{n-2}}   \sum_{i \leq q}  \a_i ^\frac{2n}{n-2}  \Big( \frac{c_7}{\l _i}\frac{\partial K}{\partial  \nu}(a_i) - c_6 \frac{\D K(a_i)}{\l_i^2}\Big) + \sum_{i =1} ^{p+q} O \Big(   \frac{1}{\l _i^3}+  \frac{ | \n K(a_i) |^2}{\l_i^2} \Big) \Big)\end{align*}
where
$$S_n := c_0^\frac{2n}{n-2} \int_{\R^n_+} \frac{ dx }{(1+| x| ^2)^n} \, ;  \, c_6:= \frac{n-2}{n^2} c_0^\frac{2n}{n-2} \int_{\R^n_+} \frac{| x |^2 dx }{(1+| x|^2)^n}   \, ;  \, c_7:= 2 \frac{n-2}{n} c_0^\frac{2n}{n-2} \int_{\R^n_+} \frac{x_n dx }{(1+| x|^2)^n} $$
\end{pro}
\begin{pf} From the definition of $J_K$, we need to expand (using the fact that $\ov{v} \perp \varphi_i$ for each $i$)
$$ \| u\|^2 = \sum \a_i^2 \| \varphi_i \|^2 + \| \ov{v} \|^2 + O\big( \sum \e_{ij}\big) = S_n \big(\sum_{i\leq q }\a _i^2+2 \sum_{i > q}\a _i^2\big) + \| \ov{v} \|^2 + O\big( \sum \e_{ij} + \sum_{ i > q} \frac{1}{\l_i^{n-2}} \big), $$
$$\int_{\mathbb{S}^n_+} K u^\frac{2n}{n-2} = \sum_{i=1}^{q+p} \a_i ^\frac{2n}{n-2} \int_{\mathbb{S}^n_+} K \varphi_i ^\frac{2n}{n-2} + \frac{2n}{n-2} \int_{\mathbb{S}^n_+} K \big( \sum \a_i \varphi_i\big)^\frac{n+2}{n-2} \ov{v} + O\Big( \sum_{i\neq j} \int \varphi_i ^\frac{n+2}{n-2} \varphi_j + \| \ov{v} \|^2 \Big).$$
The last integral is equal to $O(\e_{ij})$. The second one is presented in \eqref{777}. Concerning the first one, for $i > q$, using Lemma \ref{lem:1}, we get
\begin{align*}
\int_{\mathbb{S}^n_+} K \varphi_i ^\frac{2n}{n-2}  & = \int_{\mathbb{R}^n_+} \wtilde{K} \d_i ^\frac{2n}{n-2} + O\Big( \frac{1}{\l_i ^{(n-2)/2}}\int \d_i ^\frac{n+2}{n-2}\Big) \\
& =  \int_{B(a_i,d_i)} \wtilde{K} \d_i ^\frac{2n}{n-2} + O\Big( \int _{\R^n \setminus B(a_i, d_i) } \d_i^{\frac{2n}{n-2}} + \frac{1}{\l_i ^{(n-2)/2}}\int \d_i ^\frac{n+2}{n-2}\Big)\\
& = 2 S_n \wtilde{K}(a_i) + \frac{1}{2n} \frac{\D \wtilde{K}(a_i) }{\l_i^2} c_0^\frac{2n}{n-2} \int_{\R^n} \frac{| x | ^2 }{(1+| x | ^2)^n}dx  + O\Big( \frac{1}{\l_i^3}\Big). \end{align*}
However, for $i\leq q$, we have $\varphi_i = \d_i$ and therefore
\begin{align*}
\int_{\mathbb{S}^n_+} K \d_i ^\frac{2n}{n-2}  & = \int_{\mathbb{R}^n_+} \wtilde{K} \d_i ^\frac{2n}{n-2}
 =  S_n \wtilde{K}(a_i) + \n \wtilde{K}(a_i)   c_0^\frac{2n}{n-2} \int_{\R^n_+ } \frac{ \l_i^n (x -a_i) }{(1+\l_i^2| x-a_i | ^2)^n}dx \\
 & + \frac{1}{2} \sum \frac{\partial^2 \wtilde{K}(a_i)}{\partial x_k\partial x_\ell } c_0^\frac{2n}{n-2} \int_{\R^n_+} \frac{\l_i^n ( x -a_i)_k (x-a_i)_\ell }{(1+ \l_i ^2| x -a_i  | ^2)^n}dx  + O\Big( \frac{1}{\l_i^3}\Big)\\
 & = S_n \wtilde{K}(a_i) - \frac{\partial \wtilde{K}}{\partial \nu }(a_i) c_0^\frac{2n}{n-2} \int_{\R^n_+} \frac{x_n }{(1+| x | ^2)^n}dx +  \frac{1}{2n} \frac{\D \wtilde{K}(a_i) }{\l_i^2} c_0^\frac{2n}{n-2} \int_{\R^n_+} \frac{| x | ^2 }{(1+| x | ^2)^n}dx +  O\Big( \frac{1}{\l_i^3}\Big).
  \end{align*}
Note that, since $u\in \Sigma $, we deduce that $$J_K(u)= \frac{1}{\G^\frac{n-2}{n}} \Big( 1 + O\Big(  \sum_{i > q }  \frac{1}{\l _i^2} +\sum_{i \leq q}   \frac{1}{\l _i} | \frac{\partial K}{\partial  \nu}(a_i) | + \frac{1}{\l_i^2} +  \| \ov{v} \|^2\Big)\Big) \mbox{ where } \G := \sum_{i\leq q }\a _i^\frac{2n}{n-2}K(a_i)+2\sum_{i > q } \a _i^\frac{2n}{n-2}K(a_i).$$
Now, the precise expansion of $J_K$ follows from the above estimates, the estimate of $\|\ov{v}\|$ (see Lemma \ref{eovv}) and the fact that $(1+x)^{-(n-2)/n}= 1 -((n-2)/n) x + O(x^2)$.
\end{pf}

In the following, we will present the expansion of the gradient of $J_K$ in the potential sets. We will present the results for $p+q \geq 2$. However, the results are true for $p+q =1$, it suffices to remove the terms $\e_{ij}$'s which correspond to the interaction terms of the bubbles.
\begin{pro}\label{p:33}
Let $n\geq 5$, for $ u=\sum_{i\leq q} \a_i \d_i + \sum_{i > q} \a _j \varphi _j \in V(m,q,p,\e )$ and $i \leq q$, it holds
\begin{align*}
 \langle \n J_K(u),\l _i \frac{\partial\d _i }{\partial\l _i} \rangle  = & 2J_K(u)\Big[ - \frac{c_2} {2}\sum _{j\neq i; j\leq q} \a _j \l _i\frac{\partial \e _{ij}}{\partial \l _i}(1+o(1)) \\
 & + 2J_K(u)^\frac{n}{n-2}\a _i^\frac{n+2}{n-2}\Big( - \frac{c_3}{\l _i}\frac{\partial K}{\partial \nu}(a_i) + c_9\frac{\D K(a_i)}{\l_i ^2}\Big)\Big] +O\Big(\frac{1}{\l_i^3} + \sum_{j > q} \e_{ij} + R_1^b \Big)
\end{align*} where
 $$ R_1^b := \sum_{k \leq q} \Big(\frac{ | \n K(a_k) | }{\l _k }\Big)^{\frac{n}{2}} + \Big( \frac{1}{\l_k^2}\Big)^{\frac{n+1}{3}} + \sum_{j\neq k; j,k \leq q}  \e _{kj }^{\frac{n}{n-2}}\ln (\e _{kj }^{-1}) \, \, \,  ; \, \, \,   c_3=\frac{n-2}{2}c_0^{\frac{2n}{n-2}}\int_{\R_+^n}\frac{x_n(|x|^2-1)}{(1+|x|^2)^{n+1}}dx .$$
\end{pro}
\begin{pf}
$$ \langle \n J_K(u),\l _i \frac{\partial\d _i }{\partial\l _i} \rangle  =  2J_K(u)\Big(\sum_{j\leq q} \a_j \langle \d_j ,\l _i \frac{\partial\d _i }{\partial\l _i} \rangle  - J(u)^{n/(n-2)} \int K \big( \sum_{j \leq q} \a_j \d_j\big)^\frac{n+2}{n-2}  \l _i \frac{\partial\d _i }{\partial\l _i} + \sum_{j > q} O(\e_{ij})\Big).$$
For $j \leq q$, we have $a_j \in \partial \mathbb{S}_+^n$ and therefore, using \cite{B1}, we get, for $j\neq i$,
$$ \langle \d_j ,\l _i \frac{\partial\d _i }{\partial\l _i} \rangle = \int_{\R^n_+} \d_j^\frac{n+2}{n-2}  \l _i \frac{\partial\d _i }{\partial\l _i} = \frac{1}{2}\int_{\R^n} \d_j^\frac{n+2}{n-2}  \l _i \frac{\partial\d _i }{\partial\l _i}= \frac{1}{2} c_2 \e_{ij} + O\big(  \e _{ij }^{\frac{n}{n-2}}\ln (\e _{ij }^{-1})\big)
$$
$$ \langle \d_i ,\l _i \frac{\partial\d _i }{\partial\l _i} \rangle = \int_{\R^n_+} \d_i^\frac{n+2}{n-2}  \l _i \frac{\partial\d _i }{\partial\l _i} = \frac{1}{2}\int_{\R^n} \d_i^\frac{n+2}{n-2}  \l _i \frac{\partial\d _i }{\partial\l _i}=0.$$
Concerning the other term, it holds
\begin{align*}
\int K \big( \sum_{j \leq q} \a_j \d_j\big)^\frac{n+2}{n-2}  \l _i \frac{\partial\d _i }{\partial\l _i} =&  \sum_{j \leq q} \int K \big( \a_j \d_j\big)^\frac{n+2}{n-2}  \l _i \frac{\partial\d _i }{\partial\l _i} \\
& + \frac{n+2}{n-2} \int K  (\a_i \d_i)^\frac{4}{n-2} \big(\sum_{j \leq q; j\neq i } \a_j \d_j\big) \l _i \frac{\partial\d _i }{\partial\l _i} + O\Big( \sum_{k\neq r} \int (\d_k \d_r ) ^\frac{n}{n-2} \Big).\end{align*}
Observe that, for $j\neq i$, expanding $K$ around $a_j$, we get
\begin{align*}
 \int_{\R^n_+} K  \d_j ^\frac{n+2}{n-2}  \l _i \frac{\partial\d _i }{\partial\l _i} & = K(a_j) \int_{\R^n_+}   \d_j ^\frac{n+2}{n-2}  \l _i \frac{\partial\d _i }{\partial\l _i} +  O\Big(| \n K(a_j) |  \int_{\R^n_+} |x-a_j|   \d_j ^\frac{n+2}{n-2}  \d _i   +  \int_{\R^n_+} |x-a_j|^2   \d_j ^\frac{n+2}{n-2}  \d _i  \Big)\\
 & = K(a_j) \frac{1}{2} c_2 \e_{ij} + O\Big(  \e _{ij }^{\frac{n}{n-2}}\ln (\e _{ij }^{-1}) + \frac{| \n K(a_j) | }{\l_j}  \e _{ij }(\ln \e _{ij }^{-1}  )^{\frac{n-2}{n}} + \frac{1 }{\l_j^2}  \e _{ij }^\frac{n}{n+1} (\ln \e _{ij }^{-1}  )^{\frac{n-2}{n+1}} \Big) \\
 & = K(a_j) \frac{1}{2} c_2 \e_{ij} + O\Big(  \e _{ij }^{\frac{n}{n-2}}\ln (\e _{ij }^{-1}) + \Big(\frac{| \n K(a_j) | }{\l_j} \Big)^{n/2} + \Big( \frac{1}{\l_j^2}\Big)^{(n+1)/3}\Big)
 ,\end{align*}
\begin{align*}
 \int_{\R^n_+} K  \d_i ^\frac{n+2}{n-2}  \l _i \frac{\partial\d _i }{\partial\l _i} & =    + \sum_{k} \frac{\partial K}{\partial x_k}(a_i) \int_{\R^n_+} (x-a_i)_k \d_i ^\frac{n+2}{n-2}  \l _i \frac{\partial\d _i }{\partial\l _i} \\
 & + \frac{1}{2} \sum\frac{ \partial^2 K} {\partial x_k \partial x_\ell } (a_i)\int_{\R^n_+} (x-a_i)_k (x-a_i)_\ell   \d_i ^\frac{n+2}{n-2}  \l _i \frac{\partial\d _i }{\partial\l _i} + O\Big(  \int_{\R^n_+} |x-a_i|^3   \d_i ^\frac{2n}{n-2}  \Big)\\
 & = \frac{c_3}{\l_i}  \frac{\partial K}{\partial \nu} (a_i)  - c_9 \frac{\D K(a_i)}{\l_i^2} + O\big( \frac{1}{\l_i^3}\big).
 \end{align*}
 Finally, for $j\neq i$, it holds
 $$ \frac{n+2}{n-2} \int_{\R^n_+} K   \d_i^\frac{4}{n-2} \d_j \l _i \frac{\partial\d _i }{\partial\l _i} = K(a_i) \langle \d_j, \l _i \frac{\partial\d _i }{\partial\l _i} \rangle + O \Big( | \n K(a_i) |  \int_{\R^n_+} |x-a_i |   \d_i ^\frac{n+2}{n-2}  \d _j  + \int_{\R^n_+} |x-a_i |^2   \d_i ^\frac{n+2}{n-2}  \d _j \Big).$$
 Hence the proof follows.
\end{pf}

\begin{pro}\label{p:34}
Let $n\geq 5$. For $ u = \sum_{i\leq q} \a_i \d_i + \sum_{i > q}  \a _j \varphi _j  \in V (m ,q,p, \e )$ and $i\leq q$, it holds:
\begin{align*}
\langle \n J_K(u), \frac{1}{\l _i}\frac{\partial\d _i}{\partial a_i}\rangle &= 2J_K(u)\a _i e_n \left[ c_4\left(1- J_K(u)^{\frac{n}{n-2}}\a _i ^{\frac{4}{n-2}}K(a_i)\right)+ J_K(u)^{\frac{n}{n-2}}\a _i ^{\frac{4}{n-2}}\frac{c_5}{\l_i}\frac{\partial K}{\partial\nu}(a_i)\right]\\
&-J_K(u)c_2 \sum_{j\leq q; j\neq i} \a _j  \frac{1}{\l _i}\frac{\partial\e _{ij}}{\partial a_i}\Big(-1 + J_K(u)^{\frac{n}{n-2}}\sum_{k=i,j} \a _k ^{\frac{4}{n-2}}K(a_k)\Big)+ O\Big( \frac{1}{\l _i^{2}} \Big)\\
&-4J_K(u)^{\frac{2(n-1)}{n-2}}\a _i ^{\frac{n+2}{n-2}}\frac{2c_5}{\l _i}\n_T K(a_i) + O\Big(R_1^b + \sum_{k\leq q; k \neq i} \e _{ik}^{\frac{n+1}{n-2}}\l _k d(a_i,a_k) + \sum_{k > q } \e_{ik} \Big)
\end{align*}
where  $R_1^b$ is defined in Proposition \ref{p:33} and
$$ c_4 = (n-2)c_0^{\frac{2n}{n-2}}\int_{\R^n_+}\frac{x_n}{(1+|x|^2)^{n+1}}dx \quad \mbox{and } c_5=\frac{n-2}{2n}c_0^{\frac{2n}{n-2}}\int_{\R^n}\frac{x_n^2}{(1+|x|^2)^{n+1}}dx. $$
\end{pro}
\begin{pf}
The proof can be done as the previous one.
\end{pf}
\begin{pro}\label{p:35}
 For $ u = \sum \a _j \varphi _j \in V (m,q,p, \e )$ and $i \leq q$, we have the following expansion:
$$
\langle \n J_K(u),\d _i\rangle =2J_K(u) \a _i {S_n} \left(1-J_K(u)^{\frac{n}{n-2}}\a _i^{\frac{4}{n-2}}K(a_i)\right) + O\Big( \frac{ | \n K(a_i) | }{\l_i }+  \frac{ 1 }{\l_i^2 }+ \sum_{j\neq i} \e_{ij} \Big).
 $$
 where $S_n$ is defined in Proposition \ref{p:311}.
\end{pro}
\begin{pf}
$$ \langle \n J_K(u),\d _i\rangle =2J_K(u) \a_i \| \d_i\|^2   - J_K(u)^{n/(n-2)} \int K \d_i^\frac{2n}{n-2} +O\Big( \sum \e_{ki}\Big) . $$
Observe that
$$ \int K \d_i^\frac{2n}{n-2}  = K(a_i)  \int  \d_i^\frac{2n}{n-2} + O\Big( | \n K(a_i) |  \int | x-a_i | \d_i^\frac{2n}{n-2} + \int | x-a_i |^2 \d_i^\frac{2n}{n-2} \Big)$$ which gives the result.
\end{pf}

\begin{pro}\label{pp65}
 For $u=\sum_{j\leq q} \a_j\d_j  + \sum_{j > q} \a_j\varphi_j \in V(m,q,p , \e)$ and for each $i \geq q+1$,  we have:
\begin{align*}
& \langle\nabla J_K(u),\l_i\frac{\partial \varphi_i}{\partial\l_i}\rangle= 2J_K\biggl(-c_2 \sum_{ j\neq i}\a_j (1+o(1)) \l_i\frac{\partial\e_{ij}}{\partial\l_i} +c_2 \frac{n-2}{2} \sum_{j=q+1}^p\a_j (1+o(1)) \frac{H(a_i,a_j)}{(\l_i\l_j)^{(n-2)/2}}\\
& \qquad \qquad \qquad + c \a_i(1+o(1))\frac{\Delta K(a_i)}{\l_i^2K(a_i)}\biggr) +O\Bigl( \frac{1}{\l_i^{3}} + R_1\biggr), \\
& \langle \n J_K(u),\varphi _i\rangle =2J_K(u) \a _i {S_n}\left(1-J_K(u)^{\frac{n}{n-2}}\a _i^{\frac{4}{n-2}}K(a_i)\right) + O\Big( \frac{ | \n K(a_i) | }{\l_i }+  \frac{ 1 }{\l_i^2 }+  \frac{ 1 }{(\l_id_i)^{n-2} } + \sum_{j\neq i} \e_{ij} \Big), \\
& \langle\nabla J_K(u),\frac{1}{\l_i}\frac{\partial \varphi_i}{\partial a_i}\rangle.\nabla K(a_i)\geq c\frac{|\nabla K(a_i)|^2}{\l_i}+O\biggl(\Big(\frac{1}{\l_i^2}+ \frac{1}{(\l_id_i)^{n-2}} + \sum_{j\neq i}\e_{ij}\Big)|\nabla K(a_i)|\biggr)
\end{align*}
where $$ R_1 := \sum_{k =1}^{q+p} \Big(\frac{ | \n K(a_k) | }{\l _k }\Big)^{\frac{n}{2}} + \Big( \frac{1}{\l_k^2}\Big)^{\frac{n+1}{3}} + \sum_{j\neq k}  \e _{kj }^{\frac{n}{n-2}}\ln (\e _{kj }^{-1}) + \sum_{j > q} \frac{1}{(\l_jd_j)^n}. $$
\end{pro}

\subsection{Counting index formulae}
\begin{lem}\label{indicesz}
Let $z_1,\cdots,z_N$ be $N$ critical points of $K_1$ in $\partial \mathbb{S}^n_+$ and let $\i(z_j):= n-1 - morse(K_1,z_j)$. Assume that $$ A_1:= \sum_{j=1}^N (-1)^{i(z_j)} =1.$$
 Then the number $N$ has to be odd, say $N:= 2k +1$ (with $k\in \N_0$) and there are $k$ odd numbers $\i(z_j)$'s and $k+1$ even numbers $\i(z_j)$'s. Furthermore, for each $k\geq 0$, it hold
 \begin{align*}
 &  A_2 := \sum_{ j < \ell} (-1)^{\i(z_j) + \i(z_\ell) } = -k \quad ; \quad A_3 := \sum_{ j < \ell < r } (-1)^{\i(z_j) + \i(z_\ell) + \i(z_r) } = -k \, ,\\
 &  A_4 := \sum_{ j < \ell < r < t } (-1)^{\i(z_j) + \i(z_\ell) + \i(z_r) + \i(z_t)} = \frac{1}{2} k(k-1) .\end{align*}
\end{lem}

\begin{pf}
To compute the value of $A_2$, observe that it is the sum of $+1$ and $-1$. To get $-1$, $\i(z_j)$ and $\i(z_k)$ have to be of different parity. However, to get $+1$, $\i(z_j)$ and $\i(z_k)$ have to be of the same parity. A similar argument holds for the computation of the values $A_3$ and $A_4$. Hence:
\begin{itemize}
\item For $k=0$, we have only one point $z$ with an even $\i(z)$. Thus $A_2=A_3=A_4=0$.
\item For $k=1$, we have two points $z_0$ and $z_2$ with even $\i(z_k)$ and one point $z_1$ with an odd $\i(z_1)$. Thus, $A_4=0$, $A_3= -1$ and $A_2= 1-2=-1$.
\item For $k\geq 2$, there exist $k+1$ even numbers $\i(z_j)$ and $k$ odd numbers $\i(z_j)$. Thus,  it holds
\begin{align*}
& A_2=\begin{pmatrix} 2 \\  k+1\end{pmatrix}  + \begin{pmatrix} 2 \\ k \end{pmatrix} - \begin{pmatrix} 1 \\ k+1\end{pmatrix} \begin{pmatrix} 1 \\ k \end{pmatrix} = \frac{1}{2} (k+1)k + \frac{1}{2} k(k-1) - (k+1) k = -k,\\
 & A_3 = \begin{pmatrix} 3 \\  3\end{pmatrix} + \begin{pmatrix} 1 \\  3\end{pmatrix} \begin{pmatrix} 2 \\  2\end{pmatrix} - \begin{pmatrix} 2 \\  3\end{pmatrix} \begin{pmatrix} 1 \\  2\end{pmatrix} = -2\quad  \mbox{ if } k=2\\
  & A_3 = \begin{pmatrix} 3 \\  k+1 \end{pmatrix} + \begin{pmatrix} 1 \\  k+1\end{pmatrix} \begin{pmatrix} 2 \\  k\end{pmatrix} - \begin{pmatrix} 2 \\  k+1\end{pmatrix} \begin{pmatrix} 1 \\  k\end{pmatrix}  - \begin{pmatrix} 3 \\  k \end{pmatrix}= - k\quad  \mbox{ if } k \geq 3\\
& A_4 = \begin{pmatrix} 2 \\  3 \end{pmatrix}\begin{pmatrix} 2 \\  2\end{pmatrix}    - \begin{pmatrix} 3 \\  3 \end{pmatrix} \begin{pmatrix} 1 \\  2 \end{pmatrix}= 1 \quad  \mbox{ if } k = 2\\
   & A_4 = \begin{pmatrix} 4 \\  4\end{pmatrix} + \begin{pmatrix} 2 \\  4\end{pmatrix} \begin{pmatrix} 2 \\  3\end{pmatrix} - \begin{pmatrix} 3 \\  4\end{pmatrix} \begin{pmatrix} 1 \\  3\end{pmatrix} - \begin{pmatrix} 1 \\  4\end{pmatrix} \begin{pmatrix} 3 \\  3\end{pmatrix} = 3 \quad  \mbox{ if } k=3 \\
    & A_4 = \begin{pmatrix} 4 \\  k+1 \end{pmatrix} + \begin{pmatrix} 2 \\  k+1\end{pmatrix} \begin{pmatrix} 2 \\  k\end{pmatrix} + \begin{pmatrix} 4 \\  k\end{pmatrix} - \begin{pmatrix} 3 \\  k+1\end{pmatrix} \begin{pmatrix} 1 \\  k\end{pmatrix}  -  \begin{pmatrix} 1 \\  k+1\end{pmatrix} \begin{pmatrix} 3 \\  k \end{pmatrix}= \frac{1}{2} k (k-1) \quad  \mbox{ if } k \geq 4.
  \end{align*}
 \end{itemize}
 The proof is thereby completed.
\end{pf}

Arguing as in the above lemma, one derives the following counting formula:

\begin{lem}\label{indicesy}
Let $y_1,\cdots,y_L$ be $L$ critical points of $K$ in $ \mathbb{S}^n_+$ and let $\i(y_j):= n - morse(K,y_j)$. Assume that $$ B_1:= \sum_{j=1}^L (-1)^{\i(y_j)} =-k \quad \mbox{ with } k \geq 0.$$
 Then the number $L$ has to satisfy $L:= 2r +k$ (with $r \in \N_0$) and there are $r$ even numbers $\i(y_j)$'s and $r+k$ odd numbers $\i(y_j)$'s. Furthermore, it holds
$$   B_2 := \sum_{ 1\leq j < \ell \leq L} (-1)^{\i(y_j) + \i(y_\ell) } = -r + \frac{1}{2} k(k-1)  \, \, ; \, \,\mbox{ for each } L \geq 0.$$
\end{lem}


{\small

\bigskip
\bigskip
\bigskip
\bigskip

$$
\begin{minipage}[l]{9.5cm}
    {\small Mohameden Ahmedou}\\
    {\footnotesize Mathematisches Institut }\\
     {\footnotesize der Justus-Liebig-Universit\"at Giessen}\\
    {\footnotesize Arndtsrasse 2, D-35392 Giessen}\\
    {\footnotesize Germany}\\
    {\tt\footnotesize Mohameden.Ahmedou@math.uni-giessen.de  }
    {}
\end{minipage}
\quad
\begin{minipage}[l]{7.5cm}
    {\small Mohamed Ben Ayed}\\
    {\footnotesize Universit\'e de Sfax, }\\
      {\footnotesize Facult\'e des Sciences de Sfax}\\
    {\footnotesize D\'epartement de Math\'ematiques}\\
    {\footnotesize Route de Soukra, Sfax, BP. 1171, 3000 }\\
     {\footnotesize  Tunisia}\\
    {\tt\footnotesize  Mohamed.Benayed@fss.rnu.tn}
    {}
\end{minipage}
$$
 }

\begin{thebibliography}{99}

{ \small

\bibitem{AB}
Ahmedou, M; Ben Ayed, M. \emph{Non simple blow ups for the  Nirenberg problem on  half spheres,} preprint 2020,  arXiv:2012.11728.
\bibitem{AB20}
Ahmedou, M; Ben Ayed, M. \emph{The Nirenberg problem on  half spheres: A bubbling off analysis }, preprint 2021.


\bibitem{Au76}
 Aubin, T.  \emph{ Equations diff\'erentielles non lin\'eaires et probl\`eme de Yamabe concernant la courbure scalaire.} J. Math. Pures Appl. (9) \textbf{55} (1976), no. 3, 269–296.
\bibitem{AH91}
  Aubin, T; Hebey, E. \emph{ Courbure scalaire prescrite.} (French) [\emph{Prescribed scalar curvature}] Bull. Sci. Math. \textbf{115} (1991), no. 2, 125–131.



\bibitem{B1} Bahri A., Critical points at infinity in some variational problems, Research Notes in Mathematics, 182, Longman-Pitman, London, 1989.


\bibitem{Bahri-Invariant} A. Bahri, \emph{An invariant for yamabe-type flows with applications to scalar curvature problems in high dimensions,} A celebration of J. F. Nash Jr., Duke Math. J. \textbf{81}  (1996), 323-466.



\bibitem{BC1}
 Bahri, A;  and  Coron,J-M.  \emph{The scalar curvature problem on
the standard three dimensional spheres,} J. Funct. Anal.
\textbf{95}, (1991), 106-172.

\bibitem{BC2}
Bahri, A; Coron, J-M.  \emph{ On a nonlinear elliptic equation involving the critical Sobolev exponent: the effect of the topology of the domain.} Comm. Pure Appl. Math. \textbf{41} (1988),  253--294.



\bibitem{BLR}  Bahri, A; Li, Y.Y; Rey, O. \emph{On a variational problem with lack of compactness: the topological effect of the critical points at infinity.} Calc. Var. Partial Differential Equations \textbf{3} (1995),  67--93.

\bibitem{BN83}
 Br\'ezis, H; Nirenberg, L.   \emph{ Positive solutions of nonlinear elliptic equations involving critical Sobolev exponents.} Comm. Pure Appl. Math. \textbf{36} (1983), no. 4, 437--477.

\bibitem{BCCH}
 Ben Ayed, M; Chen,Y. ;  Chtioui, H;   Hammami,M.  \emph{On the prescribed scalar curvature problem on 4-manifolds,} Duke Mathematical Journal,
\textbf{84}, (1996), 633-677.

\bibitem{BCH}
 Ben Ayed, M; Chtioui, H; Hammami, M.  \emph{ A Morse lemma at infinity for Yamabe type problems on domains.} Ann. Inst. H. Poincar\'e Anal. Non Lin\'eaire \textbf{20} (2003),  543--577.

\bibitem{BEOA}  Ben Ayed, M; El Mehdi, K; Ahmedou, M. O. \emph{Prescribing the scalar curvature under minimal boundary conditions on the half sphere.} Adv. Nonlinear Stud. \textbf{2} (2002), no. 2, 93--116.
\bibitem{BEO}  Ben Ayed, M.; El Mehdi, K.; Ould Ahmedou, M. \emph{ The scalar curvature problem on the four dimensional half sphere,} Calc. Var. Partial Differential Equations \textbf{22} (2005), no. 4, 465--482.

\bibitem{BGO}  Ben Ayed, M; Ghoudi, R; Ould Bouh, K. \emph{Existence of conformal metrics with prescribed scalar curvature on the four dimensional half sphere.} NoDEA Nonlinear Differential Equations Appl. \textbf{19} (2012),  629--662.

\bibitem{Ben} Ben Ayed, M. \emph{ Bahri-Coron type theorem for the scalar curvature problem on high dimensional spheres}, Annali di Matematica Pura ed Applicata  \textbf{191} (2012),95–112.

\bibitem{BOA}
 Ben Ayed, M; Ould Ahmedou, M. \emph{  On the prescribed scalar curvature on 3-half spheres: multiplicity results and Morse inequalities at infinity.} Discrete Contin. Dyn. Syst. \textbf{23} (2009), no. 3, 655–683.

\bibitem{BEZ1}
Bourguignon, J.P.;  Ezin,J.P.  \emph{Scalar curvature functions in a conformal class of metrics and conformal transformations,} Trans. Amer. Math. Soc., \textbf{301},  (1987), 723-736.


\bibitem{CY}
 Chang,  A.; Yang, P.  \emph{ A perturbation result in prescribing scalar curvature on $\mathbb{S}^n$.} Duke Math. J. \textbf{64} (1991),  27--69.

\bibitem{CGY}
 Chang,  A.; Gursky, Matthew J.; Yang, Paul C. \emph{ The scalar curvature equation on 2- and 3-spheres. Calc. Var. Partial Differential Equations} \textbf{1} (1993), no. 2, 205–-229.


\bibitem{CL}  Chen, C.C; Lin, C.S. \emph{ Blowing up with infinite energy of conformal metrics on $\mathbb{S}^n.$} Comm. Partial Differential Equations \textbf{24} (1999),  785--799.

\bibitem{CL1} Chen, C.C;  Lin, C.S. \emph{Prescribing  the scalar curvature on $S^n$, I. Apriori estimates} J. Differential Geom. \textbf{57}, (2001), 67-171.

\bibitem{CL2}  Chen, C-C; Lin, C.S. \emph{Estimate of the conformal scalar curvature equation via the method of moving planes. II.} J. Differential Geom. \textbf{49} (1998), 115--178.

\bibitem{Chen-Xu}
 Chen, X; Xu, X.  \emph{The scalar curvature flow on $S^n$—perturbation theorem revisited, } Invent. Math. \textbf{187} (2012), no. 2, 395–506



\bibitem{DMOA}  Djadli, Z; Malchiodi, A; Ould Ahmedou, M. \emph{ Prescribing scalar and boundary mean curvature on the three dimensional half sphere.} J. Geom. Anal. \textbf{13} (2003), 255--289.


\bibitem{Dold}
 Dold, A.  Lectures on algebraic topology. Reprint of the 1972 edition. Classics in Mathematics. Springer-Verlag, Berlin, 1995. xii+377 pp.

\bibitem{Han}
 Han, Z-C.  \emph{ Prescribing Gaussian curvature on $\mathbb{S}^2$. } Duke Math. J. \textbf{61} (1990),  679--703.


\bibitem{KW1} J. Kazdan and F. Warner, {\it Existence and conformal deformation of metrics with prescribed Gaussian and scalar curvatures}, Ann. of Math (2) \textbf{101} (1975), 317--331.


\bibitem{LL}
Li, P.L; Liu, J.Q. \emph{ Nirenberg's problem on the two-dimensional hemi-sphere,} Int. J. Math. \textbf{4}(1993), 927-939,

\bibitem{yyli} Li, Y.Y. \emph{The Nirenberg problem in a domain with boundary,} Top. Meth. Nonlin. Anal., \textbf{6}1995, 309--329.
\bibitem{yyli1}  Li, Y.Y.  \emph{Prescribing scalar curvature on $S^{n}$ and related topics, Part I}, Journal of  Differential  Equations, \textbf{120} (1995), 319-410.


\bibitem{yyli2}
 Li, Y.Y.  \emph{Prescribing scalar curvature on $S^n$ and related
topics, Part II : existence and compactness,} Comm. Pure Appl.
Math. \textbf{49} (1996), 437-477.

\bibitem{Lions}
Lions, P.L. \emph{The concentration-compactness principle in the calculus of variations. The limit case. Part I.}  Rev. Mat. Iberoamericano \textbf{1}(1985), 145--201.


\bibitem{MM}
 Malchiodi, A;  Mayer, M.  \emph{ Prescribing Morse scalar curvatures: blow-up analysis},  preprint 2019.

 \bibitem{MM19}
Malchiodi, A; Mayer, M. \emph{ Prescribing Morse scalar curvatures: pinching and Morse theory,} Preprint 2019.

\bibitem{Schoen}
Schoen, R. Topics in Differential geometry, Graduate cours at  Stanford University, 1988 (http://sites.math.washington.edu/~pollack/research/Pollack-notes-Schoen1988.pdf)

\bibitem{SZ}
Schoen, R; Zhang, Dong.
\emph{ Prescribed scalar curvature on the n-sphere.}
Calc. Var. Partial Differential Equations \textbf{4} (1996),  1--25.

\bibitem{Struwe}
Struwe, M.  \emph{A global compactness result for elliptic boundary value problems invoving limiting nonlinearities.}  Math. Z. \textbf{187} (1984), 511--517.

}
\end{thebibliography}
\end{document}